
\def\TheMagstep{\magstep1}
\def\PaperSize{letter}          


 \long\def\TheAbstract{
 Let $X$ and $Y$ be smooth varieties of dimensions $n-1$ and~$n$ over
an arbitrary algebraically closed field, $f\:X\to Y$ a finite map that
is birational onto its image.  Suppose that $f$ is curvilinear; that
is, for all~$x\in X$, the Jacobian $\partial f(x)$ has rank at least
$n-2$.  For~$r\ge1$, consider the subscheme $N_r$ of $Y$ defined by the
$(r-1)$-th Fitting ideal of the $\O_Y$-module $f_*\O_X$, and set
$M_r:=f^{-1}N_r$.  In this setting---in fact, in a more general
setting---we prove the following statements, which show that $M_r$
and $N_r$ behave like reasonable schemes of source and target $r$-fold
points of $f$.
 \par\noindent\quad
If each component of $M_r$, or equivalently of $N_r$, has
the minimal possible dimension $n-r$,
then $M_r$ and $N_r$ are Cohen--Macaulay, and their
fundamental cycles satisfy the relation, $f_*[M_r]=r[N_r]$.  Now,
suppose that each component of $M_s$, or of $N_s$, has dimension $n-s$
for $s=1,\dots,r+1$.  Then the blowup $\Bl(N_r,N_{r+1})$ is equal to
the Hilbert scheme $\Hilb^r_f$, and the blowup $\Bl(M_r,M_{r+1})$ is
equal to the universal subscheme $\Univ^r_f$ of $\Hilb^r_f\times_Y X$;
moreover, $\Hilb^r_f$ and $\Univ^r_f$ are Gorenstein.  In addition, the
structure map $h\:\Hilb^r_f\to Y$ is finite and birational onto its
image; and its conductor is equal to the ideal~${\cal J}_r$ of
$N_{r+1}$ in $N_r$, and is locally self-linked.  Reciprocally,
$h_*\O_{\Hilb^r_f}$ is equal to $\Hom({\cal J}_r,\O_{N_{r}})$.
Moreover, $h_*[h^{-1}N_{r+1}]=(r+1)[N_{r+1}]$.\vadjust{\kern.5pt}  Similar
assertions hold for the structure map $h_1\:\Univ^r_f\to X$ if $r\ge2$.
 \par }
\def\TRUE{TRUE}
\ifx\DoublepageOutput\TRUE \def\TheMagstep{\magstep0} \fi
\mag=\TheMagstep
\parskip=0pt plus 1.75pt
\parindent10pt
\hsize26pc
\vsize42pc
\abovedisplayskip6pt plus6pt minus0.25pt
\belowdisplayskip6pt plus6pt minus0.25pt

\def\centertext
 {\hoffset=\pgwidth \advance\hoffset-\hsize
  \advance\hoffset-2truein \divide\hoffset by 2
  \voffset=\pgheight \advance\voffset-\vsize
  \advance\voffset-2truein \divide\voffset by 2
 }
\newdimen\pgwidth	\newdimen\pgheight
\def\letter{letter}	
\ifx\PaperSize\letter
	\pgwidth=8.5truein \pgheight=11truein \centertext \fi
\ifx\PaperSize\AFour
	\pgwidth=210truemm \pgheight=297truemm \centertext \fi

 \newdimen\fullhsize 	\newbox\leftcolumn
 \def\fulline{\hbox to \fullhsize}
\def\doublepageoutput
{\let\lr=L
 \output={\if L\lr
          \global\setbox\leftcolumn=\columnbox \global\let\lr=R%
        \else \doubleformat \global\let\lr=L\fi
        \ifnum\outputpenalty>-20000 \else\dosupereject\fi}%
 \def\doubleformat{\shipout\vbox{%
        \fulline{\hfil\box\leftcolumn\hfil\columnbox\hfil}%
	}%
 }%
 \def\columnbox{\vbox
   {\makeheadline\pagebody\makefootline
   \advancepageno}}%
 \fullhsize=\pgheight \hoffset=-1truein
 \voffset=\pgwidth \advance\voffset-\vsize
  \advance\voffset-2truein \divide\voffset by 2

\null\vfill\nopagenumbers\eject\pageno=1\relax 
}
\ifx\DoublepageOutput\TRUE \doublepageoutput \fi

 \font\twelvebf=cmbx12
 \font\smc=cmcsc10

\font\tenbi=cmmi10 scaled \magstep2 
 \font\sevenbi=cmmi10 \font\fivebi=cmmi7
 \newfam\bifam  \textfont\bifam=\tenbi
 \scriptfont\bifam=\sevenbi \scriptscriptfont\bifam=\fivebi
 \mathchardef\variableomega="7121 
 \mathchardef\variablenu="7117 
 \mathchardef\variabletau="711C 

\catcode`\@=11
\def\eightpoint{\eightpointfonts
 \setbox\strutbox\hbox{\vrule height7\p@ depth2\p@ width\z@}%
 \eightpointparameters\eightpointfamilies
 \normalbaselines\rm
 }
\def\eightpointparameters{%
 \normalbaselineskip9\p@
 \abovedisplayskip9\p@ plus2.4\p@ minus6.2\p@
 \belowdisplayskip9\p@ plus2.4\p@ minus6.2\p@
 \abovedisplayshortskip\z@ plus2.4\p@
 \belowdisplayshortskip5.6\p@ plus2.4\p@ minus3.2\p@
 \skewchar\eighti='177 \skewchar\sixi='177
 \skewchar\eightsy='60 \skewchar\sixsy='60
 \hyphenchar\eighttt=-1
 }
\newfam\smcfam
\def\eightpointfonts{%
 \font\eightrm=cmr8 \font\sixrm=cmr6
 \font\eightbf=cmbx8 \font\sixbf=cmbx6
 \font\eightit=cmti8
\font\eightsmc=cmcsc8
 \font\eighti=cmmi8 \font\sixi=cmmi6
 \font\eightsy=cmsy8 \font\sixsy=cmsy6
 \font\eightsl=cmsl8 \font\eighttt=cmtt8}
\def\eightpointfamilies{%
 \textfont\z@\eightrm \scriptfont\z@\sixrm  \scriptscriptfont\z@\fiverm
 \textfont\@ne\eighti \scriptfont\@ne\sixi  \scriptscriptfont\@ne\fivei
 \textfont\tw@\eightsy \scriptfont\tw@\sixsy \scriptscriptfont\tw@\fivesy
 \textfont\thr@@\tenex \scriptfont\thr@@\tenex\scriptscriptfont\thr@@\tenex
 \textfont\itfam\eightit	\def\it{\fam\itfam\eightit}%
 \textfont\slfam\eightsl	\def\sl{\fam\slfam\eightsl}%
 \textfont\ttfam\eighttt	\def\tt{\fam\ttfam\eighttt}%
 \textfont\smcfam\eightsmc	\def\smc{\fam\smcfam\eightsmc}%
 \textfont\bffam\eightbf \scriptfont\bffam\sixbf
   \scriptscriptfont\bffam\fivebf	\def\bf{\fam\bffam\eightbf}%
 \def\rm{\fam0\eightrm}%
 }

\def\vfootnote#1{\insert\footins\bgroup
 \eightpoint 
 \interlinepenalty\interfootnotelinepenalty
  \splittopskip\ht\strutbox 
  \splitmaxdepth\dp\strutbox \floatingpenalty\@MM
  \leftskip\z@skip \rightskip\z@skip \spaceskip\z@skip \xspaceskip\z@skip
  \textindent{#1}\footstrut\futurelet\next\fo@t}

\def\item#1 {\par\indent\indent\indent \hangindent3\parindent
 \llap{\rm (#1)\enspace}\ignorespaces}
 \def\part#1 {\par{\rm (#1)\enspace}\ignorespaces}

\newif\ifdates	\newif\ifdateonpageone
\def\today{\ifcase\month\or
 January\or February\or March\or April\or May\or June\or
 July\or August\or September\or October\or November\or December\fi
 \space\number\day, \number\year}
\def\today{}
\def\today{alg-geom/9412010} 

\nopagenumbers
\headline={%
 \eightpoint
  \ifnum\pageno=1\firstheadline
  \else
    \ifodd\pageno\oddheadline
    \else\evenheadline\fi
  \fi
}
 \def\firstheadline{\hfill}
 \def\oddheadline{\rm
  \rlap{January 7, 1996}
  \hfil\title\hfil\llap{\the\pageno}}
 \def\evenheadline{\rm\rlap{\number\pageno}\hfil
 	\author\hfil
	\ifdates\llap{\today}\fi}
\def\author{S. Kleiman, J. Lipman, and B. Ulrich}
\def\title{The multiple-point schemes}

\def\(#1){{\rm(#1)}}\let\leftp=(
\def\activeleftp{\catcode`\(=\active}
{\activeleftp\gdef({\ifmmode\let\next=\leftp \else\let\next=\(\fi\next}}

\def\artkey #1 {{\bf (\number\sectno.#1)}\enspace}
\def\art #1 #2.{\medbreak\noindent{\bf
 (\number\sectno.#1)\enspace}{\bf #2.}\enspace\ignorespaces}
\def\rem {\medbreak {\bf Remark }\artkey} \let\rmk=\rem
\def\dfn {\medbreak {\bf Definition }\artkey}
\def\proclaim#1 #2 {\medbreak{\bf#1 }\artkey#2 \bgroup\it\activeleftp}
\def\endproclaim{\egroup\medbreak}
\def\pf{\endproclaim{\bf Proof.}\enspace}
 \let\prp=\prop
\def\lem{\proclaim Lemma }
\def\thm{\proclaim Theorem }
\def\cor{\proclaim Corollary }

 \newcount\sectno \sectno=0
 \newskip\sectskipamount \sectskipamount=0pt plus30pt
 \def\newsect #1\par{\displayno=0 
   \advance\sectno by 1
   \vskip\sectskipamount\penalty-250\vskip-\sectskipamount
   \bigskip	
   \centerline{\bf \number\sectno. #1}\nobreak
   \medskip	
   \message{#1 }%
}

\def\dno#1${\eqno\hbox{\rm(\number\sectno.#1)}$}
\def\Cs#1){\unskip~{\rm(\number\sectno.#1)}}
\newcount\displayno
\def\eqlt#1${\global\advance\displayno by 1
 \expandafter\xdef
  \csname \the\sectno#1\endcsname{\the\displayno}
 \eqno(\the\sectno.\the\displayno)$}
\def\disp#1#2{(#1.\csname#1#2\endcsname)}
\def\Cn#1){\unskip~{\rm(\number\sectno.\csname\the\sectno#1\endcsname)}}

\newcount\MinNo
\newcount\scratch
\def\Mn#1{\scratch=\MinNo\advance\scratch by #1\relax\number\scratch}
\def\Cn{\Cs\Mn} \def\Dno{\dno\Mn}

 \newif\ifproofing \proofingfalse 
 \newcount\refno	 \refno=0
 \def\MakeKey{\advance\refno by 1 \expandafter\xdef
 	\csname\TheKey\endcsname{{%
	\ifproofing\TheKey\else\number\refno\fi}}\NextKey}
 \def\NextKey#1 {\def\TheKey{#1}\ifx\TheKey\NoKey\let\next\relax
  \else\let\next\MakeKey \fi \next}
 \def\NoKey{*!*}
 \def\RefKeys #1\endRefKeys{\expandafter\NextKey #1 *!* }

\def\UThin{\penalty\@M \thinspace\ignorespaces}
\def\relaxnext@{\let\next\relax}
\def\cite#1{\relaxnext@
 \def\nextiii@##1,##2\end@{\unskip\space{\rm[\SetKey{##1},\let~=\UThin##2]}}%
 \in@,{#1}\ifin@\def\next{\nextiii@#1\end@}\else
 \def\next{\unskip\space{\rm[{\SetKey{#1}}]}}\fi\next}
\newif\ifin@
\def\in@#1#2{\def\in@@##1#1##2##3\in@@
 {\ifx\in@##2\in@false\else\in@true\fi}%
 \in@@#2#1\in@\in@@} \def\SetKey#1{{\bf\csname#1\endcsname}}

\let\texttilde=\~
\def\~{\ifmmode\let\next=\widetilde\else\let\next=\texttilde\fi\next}

\def\p.{\unskip\space p.\penalty\@M \thinspace}
\def\pp.{\unskip\space pp.\penalty\@M \thinspace}

\catcode`\@=13

 \def\SetRef#1 #2\par{%
   \hang\llap{[\csname#1\endcsname]\enspace}%
   \ignorespaces#2\unskip.\endgraf}
 \newbox\keybox \setbox\keybox=\hbox{[18]\enspace}
 \newdimen\keyindent \keyindent=\wd\keybox
 \def\references{\vskip-\smallskipamount
  \bgroup   \eightpoint   \frenchspacing
   \parindent=\keyindent  \parskip=\smallskipamount
   \everypar={\SetRef}}
 \def\endreferences{\egroup}
 \def\nocomma{\def\eatcomma##1,{}\expandafter\eatcomma}

  \def\paper{\unskip, \bgroup\it}
 \def\paperinfo{\unskip, \rm}
  \def\inbook#1\bookinfo#2\publ#3\yr#4\pages#5
   {\unskip, \egroup in ``#1\unskip,'' #2\unskip, #3\unskip,
   #4\unskip, pp.~#5}
 \def\at.{.\spacefactor3000}
 \def\book{\unskip, \rm ``}
 \def\bookinfo#1{\unskip," #1}

 \def\nice{\inbook Enumerative and Classical Algebraic
Geometry\bookinfo P. le Barz, Y.  Hervier (eds.), Proc. Conf., Nice
1981. Progr. Math. {\bf 24} \publ Birkh\"auser \yr1982\pages}
 \def\bucharest{\inbook Algebraic Geometry, Bucharest 1982\at.
\bookinfo L. B\v adescu and D. Popescu
(eds.), Lecture Notes in Math. {\bf 1056} \publ Springer-Verlag
\yr 1984 \pages}
 \def\vancouver{\inbook Proc.~1984 Vancouver Conf.~in Algebraic
Geometry\bookinfo J.  Carrell, A. V. Geramita, P. Russell (eds.), CMS
Conf.~Proc.~{\bf 6}\publ Amer. Math. Soc. \yr 1986 \pages }
 \def\patzcuaro{\inbook Algebraic geometry and complex analysis
\bookinfo E. Ram\'\i rez de Arellano (ed.), Proc. Conf., P\'atzcuaro
1987, Lecture Notes in Math. {\bf 1414} \publ Springer-Verlag \yr1989
\pages}
 \def\sitgesii{\inbook Enumerative Geometry \bookinfo S. Xamb\'o
Descamps (ed.), Proc. Conf., Sitges 1987, Lecture Notes in Math. {\bf
1436} \publ Springer-Verlag \yr1990 \pages}

 \def\serial#1#2{\expandafter\def\csname#1\endcsname ##1 ##2 ##3
  {\unskip, \egroup #2 {\bf##1} (##2), ##3}}
 \serial{acta}{Acta Math.}
 \serial{ajm}{Amer. J. Math.}
 \serial{CR}{C. R. Acad. Sci. Paris}
 \serial{invent}{Invent. Math.}
 \serial{lmslns}{London Math. Soc. Lecture Note Series}
 \serial{mathz}{Math. Z.}
 \serial{ja}{J. Algebra}
 \serial{jlms}{J. London Math. Soc.}
 \serial{pams}{Proc. Amer. Math. Soc.}
 \serial{tams}{Trans. Amer. Math. Soc.}
 \serial{top}{Topology}

\hyphenation{eve-ry-where}

\let\:=\colon
\let\To=\longrightarrow \def\TO#1{\buildrel#1\over\To}
\def\onto{\to\mathrel{\mkern-15mu}\to}
\def\Onto{\mathrel{\setbox0=\hbox{$\longrightarrow$}%
	\hbox to\wd0{$\relbar\hss\onto$}}}
\def\smashedlongrightarrow{\setbox0=\hbox{$\longrightarrow$}\ht0=1pt\box0}
\def\risom{\buildrel\sim\over{\smashedlongrightarrow}}
\def\IP{{\bf P}} 
\def\sym{{\cal S}\!\hbox{\it ym}}

\def\Hom{{\cal H}\hbox{\it om}} \def\hom{\hbox{\rm Hom}}
\def\Im{\mathop{{\cal I}\hbox{\it m}}}
\def\Ht{\mathop{{\rm ht}}\nolimits}
\let\?=\overline \def\R#1{{\?R\,}^{#1}} \def\Rn{\R n}
\let\ox=\otimes
 
\def\Ser#1{{\rm(S${}_{#1}$)}}

\def\A{{\bf A}}
 \def\Ann{\hbox{\rm Ann}}
\def\ANN{{\cal A}\hbox{\it nn\/}}
\def\and{\hbox{ \rm and }} \def\for{\hbox{ \rm for }}
\def\Fit#1{{\cal F}\!\hbox{\it itt\/}_{#1}^}
\def\fit#1{\hbox{\it Fitt\/}^{#1}_}

\let\D=\Delta 
\def\O{{\cal O}} \def\C{{\cal C}} \def\M{{\cal M}}
\def\Idiag{{\cal I}(\Delta)}
\mathchardef\widesttilde="0367
\def\m{{\bf m}} \def\n{{\bf n}} \def\p{{\bf p}} \def\q{{\bf q}}

\def\mathopdef#1{\expandafter\def
 \csname#1\endcsname{\mathop{\rm #1}\nolimits}}
\def\NoOp{*!*}
\def\NextOp#1 {\def\TheOp{#1}\ifx\TheOp\NoOp\let\next\relax
  \else\mathopdef{#1}\let\next\NextOp \fi \next}
\NextOp
 Im Cok Ann Spec cod Sing Bl Hilb Univ H U Tor Ext rk depth codim grade
depth Supp Proj adj inf
 *!*
\def\U{{\rm U}}
\def\H{{\rm H}}
\catcode`\@=11

 \def\activeat#1{\csname @#1\endcsname}
 \def\def@#1{\expandafter\def\csname @#1\endcsname}
 {\catcode`\@=\active \gdef@{\activeat}}

\let\ssize\scriptstyle
\newdimen\ex@	\ex@.2326ex

 \def\requalfill{\cleaders\hbox{$\mkern-2mu\mathord=\mkern-2mu$}\hfill
  \mkern-6mu\mathord=$}
 \def\eqfill{$\m@th\mathord=\mkern-6mu\requalfill}
 \def\deffill{\hbox{$:=$}$\m@th\mkern-6mu\requalfill}
 \def\fiberbox{\hbox{$\vcenter{\hrule\hbox{\vrule\kern1ex
     \vbox{\kern1.2ex}\vrule}\hrule}$}}

 \font\arrfont=line10
 \def\Swarrow{\vcenter{\hbox{$\swarrow$\kern-.26ex
    \raise1.5ex\hbox{\arrfont\char'000}}}}

 \newdimen\arrwd
 \newdimen\minCDarrwd \minCDarrwd=2.5pc
 	\def\minCDarrowwidth#1{\minCDarrwd=#1}
 \def\findarrwd#1#2#3{\arrwd=#3%
  \setbox\z@\hbox{$\ssize\;{#1}\;\;$}%
  \setbox\@ne\hbox{$\ssize\;{#2}\;\;$}%
  \ifdim\wd\z@>\arrwd \arrwd=\wd\z@\fi
  \ifdim\wd\@ne>\arrwd \arrwd=\wd\@ne\fi}
 \newdimen\arrowsp\arrowsp=0.375em
 \def\findCDarrwd#1#2{\findarrwd{#1}{#2}{\minCDarrwd}
    \advance\arrwd by 2\arrowsp}
 \newdimen\minarrwd 
 \setbox\z@\hbox{$\longrightarrow$} \minarrwd=\wd\z@

 \def\harrow#1#2#3#4{{\minarrwd=#1\minarrwd%
   \findarrwd{#2}{#3}{\minarrwd}\kern\arrowsp
    \mathrel{\mathop{\hbox to\arrwd{#4}}\limits^{#2}_{#3}}\kern\arrowsp}}
 \def@]#1>#2>#3>{\harrow{#1}{#2}{#3}\rightarrowfill}
 \def@>#1>#2>{\harrow1{#1}{#2}\rightarrowfill}
 \def@<#1<#2<{\harrow1{#1}{#2}\leftarrowfill}
 \def@={\harrow1{}{}\eqfill}
 \def@:#1={\harrow1{}{}\deffill}
 \def@ N#1N#2N{\vCDarrow{#1}{#2}\UpDownarrow}
 \def\UpDownarrow{\uparrow\,\Big\downarrow}

 \def@.{\ifodd\row\relax\harrow1{}{}\hfill
   \else\vCDarrow{}{}.\fi}
 \def@|{\vCDarrow{}{}\Vert}
 \def@ V#1V#2V{\vCDarrow{#1}{#2}\downarrow}
 \def@ A#1A#2A{\vCDarrow{#1}{#2}\uparrow}
 \def@(#1){\arrwd=\csname col\the\col\endcsname\relax
   \hbox to 0pt{\hbox to \arrwd{\hss$\vcenter{\hbox{$#1$}}$\hss}\hss}}

 \def\squash#1{\setbox\z@=\hbox{$#1$}\finsm@@sh}
\def\finsm@@sh{\ifnum\row>1\ht\z@\z@\fi \dp\z@\z@ \box\z@}

 \newcount\row \newcount\col \newcount\numcol \newcount\arrspan
 \newdimen\vrtxhalfwd  \newbox\tempbox

 \def\innernewdimen{\alloc@1\dimen\dimendef\insc@unt}
 \def\measureinit{\col=1\vrtxhalfwd=0pt\arrspan=1\arrwd=0pt
   \setbox\tempbox=\hbox\bgroup$}
 \def\setinit{\col=1\hbox\bgroup$\ifodd\row
   \kern\csname col1\endcsname
   \kern-\csname row\the\row col1\endcsname\fi}
 \def\findvrtxhalfsum{$\egroup
  \expandafter\innernewdimen\csname row\the\row col\the\col\endcsname
  \global\csname row\the\row col\the\col\endcsname=\vrtxhalfwd
  \vrtxhalfwd=0.5\wd\tempbox
  \global\advance\csname row\the\row col\the\col\endcsname by \vrtxhalfwd
  \advance\arrwd by \csname row\the\row col\the\col\endcsname
  \divide\arrwd by \arrspan
  \loop\ifnum\col>\numcol \numcol=\col%
     \expandafter\innernewdimen \csname col\the\col\endcsname
     \global\csname col\the\col\endcsname=\arrwd
   \else \ifdim\arrwd >\csname col\the\col\endcsname
      \global\csname col\the\col\endcsname=\arrwd\fi\fi
   \advance\arrspan by -1 %
   \ifnum\arrspan>0 \repeat}
 \def\setCDarrow#1#2#3#4{\advance\col by 1 \arrspan=#1
    \arrwd= -\csname row\the\row col\the\col\endcsname\relax
    \loop\advance\arrwd by \csname col\the\col\endcsname
     \ifnum\arrspan>1 \advance\col by 1 \advance\arrspan by -1%
     \repeat
    \squash{\mathop{
     \hbox to\arrwd{\kern\arrowsp#4\kern\arrowsp}}\limits^{#2}_{#3}}}
 \def\measureCDarrow#1#2#3#4{\findvrtxhalfsum\advance\col by 1%
   \arrspan=#1\findCDarrwd{#2}{#3}%
    \setbox\tempbox=\hbox\bgroup$}
 \def\vCDarrow#1#2#3{\kern\csname col\the\col\endcsname
    \hbox to 0pt{\hss$\vcenter{\llap{$\ssize#1$}}%
     \Big#3\vcenter{\rlap{$\ssize#2$}}$\hss}\advance\col by 1}

 \def\setCD{\def\harrow{\setCDarrow}%
  \def\\{$\egroup\advance\row by 1\setinit}
  \m@th\lineskip3\ex@\lineskiplimit3\ex@ \row=1\setinit}
 \def\endsetCD{$\egroup}
 \def\measure{\bgroup
  \def\harrow{\measureCDarrow}%
  \def\\##1\\{\findvrtxhalfsum\advance\row by 2 \measureinit}%
  \row=1\numcol=0\measureinit}
 \def\endmeasure{\findvrtxhalfsum\egroup}

\newbox\CDbox \newdimen\sdim

 \newcount\savedcount
 \def\CD#1\endCD{\savedcount=\count11%
   \measure#1\endmeasure
   \vcenter{\setCD#1\endsetCD}%
   \global\count11=\savedcount}

 \catcode`\@=\active

\setbox0=\hbox{$\To$}		
\minCDarrowwidth\wd0		

\datestrue			
\dateonpageonefalse		

 \proofingfalse			
 \RefKeys
 AH SGAVI BVMZ BVSpgr B-E BR Cat Colley Eisenbud EGAI EGAIV GP HRD HSV
H1 H2 TJ Katz KActa KSitges KLU KU Kunz L65 L69 McCoy  MM Mats MP Mum Rob
PRob SV Ulrich ZN
 \endRefKeys
\long\def\topmatter{
\null\vskip24pt plus 12pt minus 12pt
\twelvebf			
\centerline{\title}
\centerline{of a finite curvilinear map of codimension one}
 \rm				

 \footnote{}{\noindent 1980 {\it Mathematics Subject Classification}
   (1985  {\it Revision}).  Primary 14C25; Secondary 14O20, 14N10.}
 \footnote{}{%
 {\it Acknowledgements.}\enspace It is a pleasure to thank Luchezar
Avramov, Winfried Bruns, David Eisenbud, Hans-Bj\o rn Foxby, and
Christian Peskine for fruitful discussions.  Avramov and Foxby
discussed Gorenstein maps; Bruns and Eisenbud discussed properties of
Fitting ideals; and Peskine explained at length his work \cite{GP}
with Gruson.}

 \vskip12pt plus6pt minus3pt
 \centerline{\smc
  Steven KLEIMAN,\footnote{$^{1}$}{%
    Supported in part by NSF grant 9106444-DMS.}
  Joseph LIPMAN,\footnote{$^{2}$}{%
   \hbox{Supported in part by NSA grant MDA904-92-3007}, and at MIT 21--30
   May 1989 by Sloan Foundation grant 88-10-1.}
  and Bernd ULRICH\footnote{$^{3}$}{%
   Supported in part by NSF grant DMS-9305832.}}
 \vskip15pt plus12pt minus12pt
{\parindent=24pt \narrower \noindent \eightpoint
 {\smc Abstract.}\enspace \TheAbstract
}}

 \topmatter

\newsect	Introduction

\art1 Overview.  Consider a finite map $f\:X\to Y$.  In the theory of
singularities of $f$, a leading role is played by the loci of source
and target $r$-fold points, $M_r=M_r(f)$ and $N_r=N_r(f)$.  They are
simple sets: $M_r$ is just the preimage $f^{-1}N_r$, and $N_r$ consists
of the (geometric) points $y$ of $Y$ whose fiber $f^{-1}y$ contains a
subscheme of degree~$r$ (or length $r$).  However, $M_r$ and $N_r$
support more refined structures, which reflect the multiplicities of
appearance of their points.  First of all, they support natural
positive cycles, whose classes are given, under suitable hypotheses, by
certain multiple-point formulas, which are polynomial expressions in
the Chern classes of $f$.  In fact, there are two different, but
related, derivations of these formulas: one is based on iteration
\cite{KActa}; the other, on the Hilbert scheme \cite{KSitges}.  In this
paper, we use the method of iteration to derive results about $M_r$
from corresponding results proved about $N_r$, to prove results about
the Hilbert scheme $\Hilb^r_f$ (of degree~$r$-subschemes of the fibers
of $f$), and to derive corresponding results about the universal
subscheme $\Univ^r_f$ of $\Hilb^r_f\times_Y X$.

Secondly, $M_r$ and $N_r$ support natural scheme structures: $N_r$ is
the closed subscheme of $Y$ defined by the Fitting ideal
$\Fit{r-1}Y(f_*\O_X)$, and $M_r$ is the closed subscheme $f^{-1}N_r$ of
$X$.  Under suitable hypotheses, which will be introduced in
Article~(1.3) below and developed in detail in Sections~2 and 3, these
subschemes have many lovely and desirable properties.  Work on this
matter was carried out by Gruson and Peskine \cite{GP} in 1981 and by
Mond and Pellikaan \cite{MP} in 1988.  (Beware: Mond and Pellikaan's
$M_r$ is our $N_r$; moreover, neither they nor Gruson and Peskine
really studied our $M_r$.)  In this paper, we aim to carry their work
further.  In Section~3, we establish some basic properties of the
schemes $M_r$ and $N_r$, and we prove a relation between their
fundamental cycles.  In Section~4, we relate $M_r$ and $N_r$ to
$\Hilb^r_f$ and $\Univ^r_f$ using some technical commutative algebra
developed in Section~5, the final section.

In Section~3, under suitable hypotheses, we prove, that $M_r$ and $N_r$
are `perfect' subschemes, and that their fundamental cycles satisfy the
basic relation,
	$$f_*[M_r]=r[N_r];\dno1.1$$
 see Theorems~(3.5) and (3.11).  Intuitively, this relation says that a
general point of $N_r$ has exactly $r$ preimages.  However, the
relation takes into account the multiplicity of the point as an
$r$-fold point.  In other words, the Fitting ideal gives the ``right''
nilpotent structure to the schemes $M_r$ and $N_r$.  Of course, off
$N_{r+1}$, the map $M_r\to N_r$ is finite and flat of degree~$r$ by a
standard property of the Fitting ideal.  The subtlety appears when some
component of $N_r$ is also a component of $N_{r+1}$.

For example, under suitable hypotheses, Proposition~(3.4) says that
$N_1$ is equal to the scheme-theoretic image of $f$, and that $N_2$ is
defined in $N_1$ by the conductor of $f_*\O_X$ into $\O_{N_1}$.
Therefore, with $r=2$, Relation~\Cs1.1) recovers the following
celebrated result (proved in various forms around 1950 by Apery, by
Gorenstein, by Kodaira, by Rosenlicht, and by Samuel): given the local
ring of a (closed) point of a curve on a smooth surface, the colength
of the conductor in the normalization is equal to twice the colength of
the conductor in the given ring; here $X$ is the normalized curve and
$Y$ is the smooth surface.  However, even with $r=2$, Relation~\Cs1.1)
is more general than that.  For instance, it is valid for the
birational projection into the plane of any reduced projective curve.
Furthermore, it yields the equation $\deg M_2 =2\deg N_2$ proved
indirectly by J.~Roberts \cite{Rob, 2nd par.\ p.~254} in the case of
the birational projection of a smooth projective variety of arbitrary
dimension onto a hypersurface.

 In Section 4, under suitable hypotheses, we study the Hilbert scheme
and the universal subscheme.  Notably, we prove Theorem~(4.2), which
asserts that the structure maps,
	$$h\:\Hilb^r_f\to Y \hbox{ and }h_1\:\Univ^r_f\to X,$$
 have a number of desirable properties;  also, $h^{-1}N_{r+1}$ is a
divisor, and
	$$h_*[h^{-1}N_{r+1}] = (r+1)[N_{r+1}],\dno1.2$$
 and similar assertions are valid for $h_1$ and $M_{r+1}$.  We also prove
Theorem~(4.3), which asserts the equations,
	$$\Hilb^r_f=\Bl(N_r,N_{r+1}) \hbox{ and }
	\Univ^r_f = \Bl(M_r,M_{r+1}).\dno1.3$$
 The first equation is obvious off $N_{r+1}$; indeed, $\Hilb^r_f$ is
equal to $N_r$ off~$N_{r+1}$, because $M_r\to N_r$ is finite and flat
of degree~$r$ there.  However, it is not obvious, a priori, even that
$h^{-1}N_{r+1}$ is a divisor.

Theorem~(4.4) asserts that the ideal ${\cal J}_r$ of $N_{r+1}$ in $N_r$
and the direct image $h_*\O_{\Hilb^r_f}$ are reciprocal fractional
ideals; that is, the pairing by multiplication,
	$${\cal J}_r\times h_*\O_{\Hilb^r_f}\To \O_{N_r},$$
 is well defined and perfect; in particular, ${\cal J}_r$ is the
conductor of $h_*\O_{\Hilb^r_f}$ in $\O_{N_r}$.  In addition, ${\cal
J}_r$ is locally a self-linked ideal of $\O_{N_r}$.  Furthermore, if
$r\ge2$, then similar assertions hold for $h_1$.  These results about
$h$ and~$h_1$ require the development of a lot of supporting
commutative algebra; some of it is developed at the end of Section~4,
and the rest, including a generalization of Huneke's theory of strong
Cohen--Macaulayness, is developed in Section~5.

Intuitively, the first equation of \Cs1.3) says that when $N_r$ is
blown up along $N_{r+1}$, then the $(r+1)$-fold points of $f$ on $Y$
are resolved into their constituent $r$-fold points.  Equation \Cs1.2)
says that the number of constituents is $r+1$, just as there should be
since there are $r+1$ different ways in which a group of $r$ points can
be chosen from a group of $r+1$ points.  Similar statements hold for
the $r$-fold points of $f$ on $X$.  Moreover, the second equation of
\Cs1.3) formally implies the equation,
	$$[M_r]=h_{1*}[\Univ^r_f],$$
 which says that $[M_r]$ is equal to the cycle whose class is given by
the refined $r$-fold-point formula of \cite{KSitges, (1.18), p.~107}.
Furthermore, the first equation of \Cs1.3) implies that the
$r$-fold-point formula is valid when $N_s$ has codimension $s$ for
$s=1,\dots,r$ (assuming the usual hypotheses on~$f$ and $Y$ in
addition).  Thus the present paper clarifies the enumerative
significance of the refined $r$-fold-point formula.

\art2 Applications.  The theory in this paper applies, for example, to
the enumeration of the secant curves of a given space curve $C$.
Indeed, Gruson and Peskine \cite{GP} made their development of the
theory to give modern derivations of the nineteenth century formulas
for the degree and genus of the curve of trisecant lines of $C$, and
for the number of quadrisecant lines.  They used the following setup:
$Y$ is the Grassmannian of lines; $X$ is the incidence variety of pairs
$(P,L)$ where $P$ is a point of $C$, and $L$ is a line through $P$; and
$f\:X\to Y$ is the projection.  Then $N_r$ is the scheme of $r$-secant
lines.  So the degree of $N_4$ is the number of quadrisecants, and it
is given by a formula of Cayley's.  To obtain this formula, Gruson and
Peskine used the Grothendieck--Riemann--Roch theorem and Porteous's
formula; however, instead, it's possible to use the 4-fold-point
formula.

Similarly, by using a stationary multiple-point formula, Colley
\cite{Colley, 5.8, p.~62} recovered Salmon's formula for the number of
reincident tangent lines of $C$.  In much the same way, S. Katz
\cite{Katz, 2.5, p.~151} recovered Severi's formula for the number of
8-secant conics to $C$: he worked out the 8-fold-point formula for the
map $f\:X\to Y$, where $Y$ is the variety of conics in space, and where
$X$ is the incidence variety of pairs $(P,L)$ where $P$ is a point of
$C$ and $L$ is a conic through $P$.  Later, Johnsen \cite{TJ}
established the enumerative significance of Severi's formula for
curves~$C$ of two types: (1)~complete intersections of general pairs of
surfaces, each of degree at least 15, and (2)~general rational curves
of degree at least~15.  He did so by using Gruson and Peskine's local
analysis as a model to show that in each case $N_r$ has codimension $r$
for all $r$; then he referred to a preliminary version of the present
paper, and quoted the discussion given at the end of Article~(1.1)
above to complete the proof.

Gruson and Peskine obtained the geometric genus of $N_3$ as follows
\cite{GP, Thm.~3.6, p.~25}: first they found its arithmetic genus; then
they determined that, under blowing-up along $N_4$, the arithmetic
genus drops by the amount of $3\deg N_4$; and finally they proved a
necessary and sufficient geometric condition for this blowup to be
smooth.  To determine the amount of the drop in the genus, they proved
an abstract algebraic result \cite{GP, 2.6, p.~13} and a related
geometric result \cite{GP, 2.7, p.~14}.  The former applied, a priori,
to a certain modification of $N_3$, and the latter result identified
this modification as the blowup of $N_3$ along $N_{4}$.  Now, this
blowup is equal to $\Hilb^3_f$ by \Cs1.3).  Correspondingly, their
algebraic result, which they proved for an arbitrary $r$, becomes
simply Equation~\Cs1.2).  Thus we recover their result.  In fact, we
derive \Cs1.2) from \Cs1.1) by induction on $r$.  Initially, the two
results coincide as $\Hilb^1_f=X$ and $h=f$; more precisely,
\Cs1.2) with $r=1$ coincides with \Cs1.1) with $r=2$.  Gruson and
Peskine \cite{GP, p.~13} themselves observed that, when $r=1$, they had
recovered the old result about the colength of the conductor of a curve
on a smooth surface.  Thus \Cs1.2) and \Cs1.1) may be viewed as
different generalizations of this old result.

It is of some importance to determine how $M_r$ and $N_r$ vary in a
family.  Of course, since they are defined by the $(r-1)$-th Fitting
ideal of the $\O_Y$-module $f_*\O_X$, their formations commute with
base change.  So the problem is to find conditions guaranteeing that
these schemes are flat when $X$ and $Y$ are flat over a given parameter
space.  Mond and Pellikaan devoted much of their paper \cite{MP} to the
issue; they considered only $N_r$, but the situation is similar for
$M_r$.  Assuming that the parameter space is smooth, they noted
\cite{MP, top p.~113} that $N_r$ is flat if it is Cohen--Macaulay and
its fibers are equidimensional of constant dimension.  Although $N_r$
is defined by a Fitting ideal, the expected codimension~$r$ of $N_r$ is
not the ``generic'' value for that Fitting ideal.  So a portion of
\cite{MP} is devoted to re-expressing the ideal of $N_r$, locally, as
the zeroth Fitting ideal of a suitable $\O_Y$-module under suitable
hypotheses; see \Cs3.1).  Then they could conclude that $N_r$ is
Cohen-Macaulay.

Similarly, we prove Theorem~(3.5) below by re-expressing the ideal
of~$N_r$ as a zeroth Fitting ideal; this theorem asserts, in particular,
that $N_r$ is perfect of grade $r$.  We derive the corresponding result
for $M_r$ from this result for $N_r$ in Theorem~(3.11).  These results
yield the flatness of $M_r$ and $N_r$, by virtue of the local
criterion, without any special assumptions on the parameter space.

\art3 Hypotheses. In this paper, we work with a finite map $f\:X\to Y$
of arbitrary locally Noetherian schemes (whereas Gruson and Peskine
\cite{GP} worked with algebraic varieties, and Mond and Pellikaan
\cite{MP} worked with complex analytic spaces, although much in their
papers generalizes with little or no change).  Thus our results
apply not only to individual varieties in arbitrary characteristic, but
also to families of varieties, including infinitesimal families and
families of mixed characteristic.  Moreover, for the most part, there
would be little technical advantage in it if we worked over a field,
let alone over an algebraically closed field or over a field of
characteristic 0.

To develop the theory fully, we need to make a number of hypotheses.
Often we need to assume, for an appropriate $r$, that $Y$ satisfies
{\it Serre's condition\/} \Ser r; that is, every local ring $\O_{Y,y}$
is either Cohen--Macaulay or of depth at least $r$.  In addition, we
make six hypotheses on $f\:X\to Y$.  However, they are not independent.
We now discuss these six, one after the other.

The first hypothesis on $f$ is that $f$ be {\it locally of
codimension~$1$}; in other words, every local ring $\O_{X,x}$ is of
dimension $1$ less than that of $\O_{Y,fx}$.  Without this hypothesis,
the $N_r$ needn't be Cohen--Macaulay when they should be.  To
illustrate this point, Mond and Pellikaan \cite{MP, bottom p.~110} gave
the following example: $X$ is the $t$-line; $Y$ is $3$-space; and
$f(t):=(t^3,t^4,t^5)$.  Here $N_1$ is not Cohen--Macaulay at the origin
because it is not reduced there.  However, this $f$ is not
dimensionally generic, as $f$ is singular at the origin; so $N_2$ has
codimension 1 in $N_1$, whereas its expected codimension is the
codimension of $f$, namely, 2.  On the other hand, Joel Roberts
(pvt.\space comm., April 18, 1991, see also \cite{ZN, Cor.~2.7})
gave the following argument, which shows that the preceding phenomenon
is not accidental.  Suppose that $X$ and $N_1$ are both
Cohen--Macaulay, and consider the sheaf
	$$\M_2:=\O_{f_*\O_X}/\O_{N_1}.$$
 Its support is the set $N_2$ and, at any point $x$ of $N_2$, its depth
is at least $\depth(\O_{X,x})-1$, which is equal to $\dim(\O_{X,x})-1$.
However, the depth of a sheaf is at most the dimension of its support.
Hence, at $x$, the codimension of $N_2$ in $N_1$ cannot be 2 or more.
However, again, its expected value is the codimension of $f$.

The second hypothesis is that $f$ be {\it locally of flat dimension
$1$}; in other words, every local ring $\O_{X,x}$ is an
$\O_{Y,fx}$-module of flat dimension $1$.  It is equivalent that
$f_*\O_X$ be presented, locally, by a square matrix with regular
determinant; see \cite{MP, 2.1, p.~114} and Lemma~(2.3) below.  By the
same token, it is equivalent that $N_1$ be a divisor.  Now, the second
hypothesis implies the first by Corollary~(2.5).  In practice, often
$Y$ is nonsingular; if so, then, by the Auslander--Buchsbaum formula,
the second hypothesis obtains if and only if the first does and $X$ is
Cohen--Macaulay.

The third hypothesis is that $f$ be {\it birational\/} (or of degree~1)
onto its image.  Suppose that $Y$ satisfies Serre's condition \Ser 2,
and that the first three hypotheses obtain.  Then Proposition~(3.4)
implies that $N_1$ is equal to the scheme-theoretic image of $f$; in
other words, the Fitting ideal $\Fit0Y(f_*\O_X)$ is equal to the
annihilator $\ANN_Y(f_*\O_X)$.  Moreover, then $N_2$ has codimension 2
in $Y$, and
  $$\ANN_Y(f_*\O_X/\O_{N_1})=\Fit1Y(f_*\O_X)
  =\Fit0Y(f_*\O_X/\O_{N_1}).$$
 The history of these equations is involved, and was indicated in the
discussion of (1.6), (1.7) and (1.8) in \cite{KLU, p.~202}; since then,
Zaare-Nahandi \cite{ZN} handled a few additional, but special, cases,
in which $f$ need not have codimension 1.  The first equation above
implies that $N_2$ is defined in $N_1$ by the conductor; the second
equation implies that $N_2$ is perfect (compare also with J.~Roberts'
\cite{Rob, Thm.~3.1, p.~258}).

The fourth hypothesis is that f be {\it curvilinear;} in other words,
the differential corank of $f$, that is, the corank of the Jacobian map
  $$\partial f(x)\: f^*\Omega ^1_{Y}(x) \To \Omega ^1_{X}(x),$$
 is at most 1 at every $x$ in $X$.  This hypothesis implies that $f$ is
{\it cyclic;} that is, locally the $\O_Y$-algebra $f_*\O_X$ has a
primitive element.  This implication was proved in the case that $X$ is
smooth by Marar and~Mond \cite{MM, 2.9, p.~560}, and it is proved in
full generality in Proposition~(2.7) below.   As a special
case, this proposition contains the usual theorem of the primitive
element for a field extension with limited inseparability.

Assume that $f$ is cyclic and that $Y$ satisfies Serre's condition \Ser
r.  If, in addition, $N_r$ locally has codimension $r$ in $Y$, then
$N_r$ is perfect by Theorem~(3.5); in fact, if $a$ is a primitive
element at $y\in N_r$, then
	$$\Fit{r-1}Y(f_*\O_X)_y =\Fit0Y (\M_r)\dno3.1$$
 where
	$$\M_r:=
    (f_*\O_X)_y\Big/\sum\nolimits_{i=0}^{r-2}\O_{Y,y}a^i.$$
 This relation between Fitting ideals was proved by Mond and Pellikaan
\cite{MP, 5.2, p.~136} under the additional assumption that $N_{r+1}$
has codimension $r+1$.  Briefly put, they proved that the two ideals
are equal off $N_{r+1}$, where the job is simpler because $(f_*\O_X)_y$
is generated by $1,\dots,a^{r-1}$; then they concluded that the two
ideals are equal everywhere because the one on the left contains the
one on the right, and the latter has no embedded components.  However,
the relation had already been proved without the assumption on the
codimension of $N_{r+1}$ by Gruson and Peskine \cite{GP, 1.3, p.~4};
they gave an elementary argument which applies to any finite cyclic
extension of an arbitrary commutative ring.  The relation plays an
essential role in the present paper, entering in the proofs of
Lemma~(3.6) and Theorem~(4.4).

In passing, let's note some other interesting properties of $\M_r$.
First,
    $$\Proj(\sym(\M_r)/\O_{N_r}\hbox{-torsion})=\Bl(N_r,N_{r+1}).$$
 if $Y$ satisfies \Ser{r+1}, if $N_r$ and $N_{r+1}$ are locally of
codimensions $r$ and $r+1$ in $Y$, and if $f$ is also locally of flat
dimension 1.  This equation was proved by Gruson and Peskine \cite{GP,
2.7, p.~14} for $r=3$, and here in brief is a version of their proof
for arbitrary $r$.  There is a natural surjection $\mu$ from $\M_r$
onto the sheaf associated to the ideal $I$ in the proof of
Theorem~(4.4) below.  In that proof, it is shown that $\Bl(I)$ is equal
to $\Bl(N_r,N_{r+1})$.  Now, $\mu$ is, obviously, an isomorphism modulo
$\O_{N_r}$-torsion; in fact, $\mu$ is an isomorphism because the
(Fitting) ideals defining $N_r$ and $N_{r+1}$ have the correct grades,
namely, $r$ and $r+1$.  Thus the equation holds.  Second, there is a
natural surjection from $\sym_{r+1}{\M_{r}}$ onto the ideal ${\cal
J}_r$ of $N_{r+1}$ in $N_r$, as Gruson and Peskine \cite{GP, 2.2, p.~8}
show, and it factors through $\sym_{r+1}(\mu)$.

The fifth hypothesis is that $f$ be {\it Gorenstein}; that is, $f$ has
finite flat dimension and its dualizing complex $f^!O_Y$ is
quasi-isomorphic to a shifted invertible sheaf.  If also the second
hypothesis obtains (that is, $f$ is locally of flat dimension 1), then
$f_*\O_X$ is presentable locally by a symmetric matrix.  This result
was proved over the complex numbers by Mond and Pellikaan \cite{MP, 2.5
p.~117}, and a version of it had been proved earlier by Catanese
\cite{Cat, 3.8, p.~84}.  Mond and Pellikaan \cite{MP, 4.3 p.~131} went
on to prove that, if $N_3$ has codimension in $Y$ at least 3, then
$N_3$ has codimension exactly 3 and $N_3$ is Cohen--Macaulay because
then $\Fit{2}Y(f_*\O_X)_y$ is locally a symmetric determinantal ideal.
Those results will not be recovered in this paper; however, see
\cite{KU} where there are new proofs, which, unlike the old, work in
the present general setting, and there are related new proofs of the
characterization (due to Valla and Ferrand) of perfect self-linked
ideals of grade 2 in an arbitrary Noetherian local ring as the ideals
of maximal minors of suitable $n$ by $n-1$ matrices having symmetric
$n-1$ by $n-1$ subblocks.

The sixth hypothesis is that $f$ be a {\it local complete
intersection}; that is, each point of $X$ has a neighborhood that is a
complete intersection in some (and so any) smooth $Y$-scheme (see
\cite{SGAVI, VIII 1,\pp.466--75}).  For example, if $X$ and $Y$ are
smooth, then the graph map of $f$ embeds $X$ in $X\times Y$; so $f$ is
a local complete intersection.  In the presence of the second and
fourth hypotheses, the fifth and sixth are equivalent. Indeed, the
sixth obviously implies the fifth; the converse is proved in
Proposition~(2.10) via a version of an old argument of Serre's.  (The
converse should now be borne in mind when reading \cite{KLU, bot.\
p.~200}.)  Moreover, by Proposition~(2.10), the sixth hypothesis,
combined with the first, implies the second.  The sixth hypothesis is
required in global enumerative multiple-point theory to ensure an
adequate theory of Chern classes of $f$ and of the pullback operator
$f^*$.  Although the sixth hypothesis played no special role in Mond
and Pellikaan's paper \cite{MP}, the hypothesis plays an essential role
in the present paper.

The various hypotheses on $f$ are inherited by the iteration map $f_1$
thanks to Lemma~(3.10), and this fact plays a leading role in this
paper.  The {\it iteration map} $f_1\:X_2\to X$ is defined as follows:
$X_2:=\IP(\Idiag)$ where $\Idiag$ is the ideal of the diagonal, and
$f_1$ is induced by the second projection.  It is remarkable how strong
a condition it is for $f_1$ to satisfy the second hypothesis; indeed,
by Proposition~(3.12), if $f_1$ does, if $Y$ satisfies \Ser2, and if
$f$ satisfies the second and fourth hypotheses, then $f$ satisfies all
six.  The usefulness of $f_1$ stems from the equations,
	$$M_r(f) = N_{r-1}(f_1) \hbox{ and }
	\Univ^r_f = \Hilb^{r-1}_{f_1}\for r\ge2,$$
 which are proved in Lemma~(3.9) and Proposition~(4.1).  These
equations permit us to derive general properties of the source
multiple-point loci and the universal subschemes from corresponding
properties of the target multiple-point loci and the Hilbert schemes,
proceeding by induction on~$r$ when convenient.

\newsect	Special finite maps

\dfn1
 A map $f\:X\to Y$ of schemes will be said to be {\it locally of flat
dimension s\/} if $X$ is nonempty and if, for every $x$ in $X$, the
local ring $\O_{X,x}$ is of flat dimension $s$ over $\O_{Y,fx}$.

\prp2
 Let $f\:X\to Y$ be a finite map that locally has a finite
presentation.  Then $f$ is locally of flat dimension~$1$ if and only if
its scheme of target points $N_1$ is a divisor of $Y$.
 \pf
 The assertion results immediately from the equivalence of (iii) and
(v) in the following  lemma, where the rings need not be Noetherian.

\lem3
 Let $\phi\:R\to B$ be a nonzero homomorphism of rings.  Assume that
$B$ has a finite presentation as an $R$-module. Then the following five
conditions are equivalent: \smallskip
 \item i The ring $B$ has flat dimension~$1$ over $R$, and for every
prime~$\p$  in $B$, the prime~$\q:=\phi^{-1}\p$ is nonminimal in $R$.
 \item ii For every prime~$\p$ in $B$, the prime~$\q:=\phi^{-1}\p$ is
nonminimal in $R$, and the localization $B_{\p}$ has flat dimension\/ {\rm
at most} $1$ over $R_{\q}$.
 \item iii For every prime~$\p$ in $B$, the localization $B_{\p}$ has flat
dimension {\rm exactly} $1$ over $R_{\q}$ where $\q:=\phi^{-1}\p$.
 \item iv For every maximal ideal $\q$ of $R$, the $R_{\q}$-module
$B_{\q}:=B\ox R_{\q}$ is presented by a square matrix whose determinant is
regular.
 \item v The zeroth Fitting ideal $\Fit0R(B)$ is invertible.
 \smallskip\noindent
 The preceding (equivalent) conditions imply the following condition,
and they are all equivalent if $R$ is Noetherian. \smallskip
 \item vi The $R$-module $B$ is perfect of grade 1.
 \smallskip\noindent
 Furthermore, if (iv) obtains, then a suitable square matrix may be
obtained from any matrix presenting~$B$ by omitting suitable columns.
 \pf
 Assume (i).  Now, for any $R$-module $M$,
	$$B_{\p}\ox_B\Tor_i^R(B,M) =  \Tor_i^{R_{\q}}(B_{\p},M_{\q});\dno3.1$$
 indeed, for any free resolution $E_\bullet$ of $M$,
	$$B_{\p}\ox_B H_i(B\ox_R E_\bullet)= H_i(B_{\p}\ox_B B\ox_R
	  E_\bullet) = H_i(B_{\p}\ox_{R_{\q}}R_{\q}\ox_R E_\bullet).$$
  Hence $B_{\p}$ has flat dimension at most 1 over $R_{\q}$, because every
$R_{\q}$-module $N$ is of the form $M_{\q}$ (for example, take $M:=N$).  Thus
(ii) holds.

Assume (ii).  Then the minimal primes $\p'$ in $B_{\p}$ have nonminimal
preimages $\q'$ in $R_{\q}$.  So, if $\q''$ is a minimal prime contained in
$\q'$ and if $a\in (\q'-\q'')$, then $a$ is regular on $R_{\q}/\q''$, but not
on $B_{\p'}/\q''B_{\p'}$, since its image in $B_{\p'}/\q''B_{\p'}$ is
nilpotent.
Hence
$B_{\p}$ is not flat over $R_{\q}$.  Therefore, $B_{\p}$ has flat dimension
exactly 1 over $R_{\q}$.  Thus (iii) holds.

Assume (iii).  To prove (iv), we may assume that $\q\supset\ker \phi$,
because otherwise $B_{\q}=0$.  Then $\q$ is of the form $\phi^{-1}\p$
because $\phi$ is finite.  So \Cs3.1) implies that $B_{\q}$ has flat
dimension at most 1 over $R_{\q}$.  By hypothesis, there is a short exact
sequence of $R_{\q}$-modules,
	$$0\To E\To F\To B_{\q}\To 0,$$
 in which $F$ is free and $E$ is finitely generated.  Hence $E$ is
flat, and therefore free by, for example, \cite{Mats, 7.10. p.~51}.

Set $I:=\Fit0R(B)$.  Then $\Ann(I_{\q})$ vanishes by McCoy's theorem
\cite{McCoy, Thm.~51, p.~159} (or by \cite{L65, Lem., p.~889}) because
$E$ is free.  Hence, $\q$ is not a minimal prime of $R$; indeed,
otherwise, $R_{\q}$ would have dimension~0, so $B_{\q}$ would be free
because $I_{\q}$ would be equal to $R_{\q}$, but $B_{\p}$ has flat
dimension exactly 1.  Hence, if $\q'$ is any minimal prime of $R$
contained in $\q$, then $B_{\q}\ox R_{\q'}=0$.  Therefore, $\rk E=\rk
F$.  In other words, $B_{\q}$ is presented by a square matrix ${\bf
M}$.  Now, $\det{\bf M}$ generates $I_{\q}$, and $\Ann(I_{\q})=(0)$; so
$\det{\bf M}$ is regular. Thus (iv) holds.

Obviously (iv) implies (v).  The converse is a special case of
\cite{L69, Lem.~1, p.~423}, but may be proved directly as follows.
Assume (v).  Take any matrix ${\bf M}$ presenting~$B_{\q}$ over $R_{\q}$; say
${\bf M}$ is $m$ by $n$ with $m\le n$.  Moreover, we may assume that
the zeroth Fitting ideal is generated by the determinant of the
submatrix ${\bf N}$ formed by the first $m$ columns because the ideal
is invertible and $R_{\q}$ is local.  Hence $\det {\bf N}$ is regular, and
divides the determinant of every $m$ by $m$ submatrix of~${\bf M}$.
Hence, by Cramer's rule, every column of ${\bf M}$ is a linear
combination of the first~$m$.  Therefore, ${\bf N}$ too presents $B_{\q}$.
Thus (iv) and the last assertion hold.\looseness=-1

Assume (iv).  Then, for every prime~$\q$ of $R$, the $R_{\q}$-module $B_{\q}$
has flat dimension at most~1.  Hence the $R$-module $B$ has flat
dimension at most~1.  If $\q=\phi^{-1}\p$ for some prime $\p$ in $B$, then
$\q R_{\q}$ contains a regular element, namely, the determinant of a matrix
presenting~$B_{\q}$; indeed, this determinant lies in $\Ann_{R_{\q}}B_{\q}$,
which is contained in~$\q R_{\q}$.  Hence, $\q$ is not minimal.  Moreover,
$\hom_{R_{\q}}(B_{\q},R_{\q})=0$; whence, $\hom_{R}(B,R)=0$ because $B$ has a
finite presentation.  Now, if $\p$ is a minimal prime of~$B$, if $\q'$ is
a minimal prime contained in $\q=\phi^{-1}\p$ and if $a\in(\q-\q')$, then
$a$ is regular on $R/\q'$, but not on $B/\q'B$; hence $B$ is not flat
over~$R$.  Alternatively, $B$ is not flat over $R$ because its zeroth Fitting
ideal is invertible, so nonzero.  Thus both (i) and (vi) hold.

Finally, assume that (vi) holds and that $R$ is Noetherian.  Then
$\Ann_RB$ contains a regular element.  This element lies in any prime~$\q$
of~$R$ of the form $\q=\phi^{-1}\p$ where $\p$ is a prime of $B$.
Hence, $\q$ is not minimal.  Thus (i) holds, and the proof is complete.

\dfn4
 Let $f\:X\to Y$ be a map of locally Noetherian schemes, and $s$ an
integer.  Call $f$ {\it locally of codimension $s$\/} if $X$ is
nonempty and if, for every $x$ in $X$,
	$$\dim\O_{X,x}=\dim\O_{Y,fx}-s.$$

\cor5
 Let $f\:X\to Y$ be a finite map of locally Noetherian schemes.  Assume
that $f$ is locally of flat dimension~$1$.
 \part1
 Then $f$ is locally of codimension~$1$.
 \part2
 Then the fundamental cycles satisfy the relation, $f_*[X]=[N_1]$.
 \pf{\def\R{\widehat R}\def\B{\widehat B}\def\C{\widehat C}%
 To prove (1), let $x$ be a point of $X$, and $y$ its image in $Y$.  Let
$R$ and $B$ be the corresponding local rings, and $C$ the semi-local
ring of $f^{-1}y$.  Let `$\widehat{\phantom{a}}$' denote completion.
It suffices to check that $\dim\B=\dim\R-1$.  Now,
$\Fit0RC$ is invertible by \Cs2); hence, $\Fit0{\R}\C$ is
invertible.  However, $\B$ is a direct summand of $\C$.
Hence $\Fit0{\widehat R}{\widehat B}$ is invertible.  Therefore,
$\dim({\widehat R}/\Fit0{\widehat R}{\widehat B})$ is equal to
$\dim{\widehat R}-1$.  Finally, $\R/\Fit0{\R}\B$ and $\B$ have
the same dimension because $\Fit0{\R}\B$ and $\Ann_{\R}(\B)$ have the
same radical.  Thus (1) holds.
  }

Consider (2).  Since $N_1$ is a divisor by \Cs2), at the generic point
$\nu$ of any of its components, the lengths $\ell_{\nu}(f_*\O_X)$ and
$\ell_{\nu}(\O_{N_1})$ are equal by the determinantal length formula,
\cite{BR, 2.4, 4.3, 4.5} or \cite{BVMZ, 2.10, p.~154}; in other
words, (2) holds.

\dfn6
 Let $f\:X\to Y$ be a map of schemes, $x\in X$.  Call the number,
	$$\dim_{k(x)}\Omega^1_f(x),$$
 the {\it differential corank\/} of $f$ at $x$.  In terms of a
base scheme $S$, this number is simply the corank of the Jacobian map,
  $$\partial f(x)\: f^*\Omega ^1_{Y/S}(x) \To \Omega ^1_{X/S}(x).$$
 Call $f$ {\it curvilinear\/} if its differential corank is at most~1
at every $x$ in $X$.

\prp7
 Let $f\:X\to Y$ be a finite map of schemes, $t\ge1$.
 \part1 Let $x\in X$.  Then $f$ has differential corank at most~$t$ at
$x$ if and only if $x$ has a neighborhood $U$ such that the restriction
$f|U$ factors through an embedding of $U$ into the affine $t$-space
$\A^t_Y$.
 \part2 Let $y\in Y$ be a point whose residue class field is infinite.
Then $f$ has differential corank at most~$t$ at every point $x$ of
$f^{-1}y$ if and only if $y$ has a neighborhood $V$ such that the
restriction $f^{-1}V\to V$ factors through a closed embedding of
$f^{-1}V$ in the affine $t$-space $\A^t_V$.
 \pf
 The assertions follow immediately from the next lemma.

\lem8
 Let $R$ be a local ring \(which need not be Noetherian).  Let $B$ be
an $R$-algebra that's finitely generated as an $R$-module, $t\ge1$.
 \part1 Let $\m$ be a maximal ideal of $B$, let $C$ be an $R$-subalgebra
of $B$ generated by $t$ elements, and set $\n:=\m\cap C$.  If the
canonical map $C_{\n}\to B_{\m}$ is surjective, then
	$$\dim_{B/\m}(\Omega^1_{B/R}/\m\Omega^1_{B/R}) \le t\dno8.1$$
 \part2 Let $\m$ be a maximal ideal of $B$ such that \Cs8.1) obtains.
Then there exists an $R$-subalgebra $C$ of $B$ generated by $t$ elements
such that, if $\n:=\m\cap C$, then the natural map $C_{\n}\to B_{\m}$ is
bijective.
 \part3 Assume that \Cs8.1) obtains for every maximal ideal $\m$ of $B$,
and that the residue class field of $R$ is infinite.  Then $B$ is
generated as an $R$-algebra by $t$ elements.
 \pf
 In (1), say $C$ is generated by $x_1,\dots,x_t$.  Then the images of
$dx_1,\dots,dx_t$ generate $(\Omega^1_{C/R})_{\n}$.  Hence \Cs8.1) holds,
as asserted.

 To prove (2) and (3), we may assume that $R$ is a field.  Indeed, let
$k$ be the residue field of $R$.  First, consider (3).  Suppose that
$B\ox k$ is generated by $t$ elements.  Lift them to $B$, and let $C$
be the resulting subalgebra.  Then the inclusion $C\to B$ is surjective
by Nakayama's lemma, because $C\ox k\to B\ox k$ is surjective by
assumption and $B$ is finitely generated as an $R$-module by
hypothesis.  Thus, to prove (3), we may replace $R$ and $B$ by $k$ and
$B\ox k$.

Next, consider (2).  Suppose that there exists a subalgebra $C$ of $B$
generated by $t$ elements, such that if $\n:=\m\cap C$, then $C_{\n}\ox
k\to B_{\m}\ox k$ is surjective.  If $x_1$ lies in $\n$, replace it by
$x_1+1$; obviously, doing so does not change $C$.  Now, $B\ox k$ is
equal to the direct product of its localizations at the maximal ideals
of $B$; so there is an element $x$ of~$B$ whose image in $B_{\m}\ox k$ is
equal to 1 and whose image in the other localizations is equal to 0.
For each $i$, replace $x_i$ by $x_ix$; doing so may change $C$, but
$C_{\n}\ox k\to B_{\m}\ox k$ remains surjective.  The following argument,
adapted from \cite{EGAIV, (18.4.6.1), p.~120}, shows that $C_{\n}\to
B_{\m}$ is now bijective.

The map $C_{\n}\to B_{\m}$ factors as follows:
	$$C_{\n}\TO\alpha B_{\n}\TO\beta B_{\m}.$$
 Consider $\beta$.  It will be bijective if every element of $B-\m$
becomes a unit in $B_{\n}$, so if every maximal ideal $\p$ of $B_{\n}$
contracts to $\m$.  Since $\alpha$ is finite, $\alpha^{-1}\p$ is a
maximal ideal of $C_{\n}$; hence, $\alpha^{-1}\p=\n C_{\n}$.
Therefore, if $\q$ denotes the trace of $\p$ in $B$, then $\q\cap
C=\n$.  Since $B$ is a finitely generated $C$-module, $\n$ is a maximal
ideal of $C$, and so $\q$ is a maximal ideal of $B$.  Now, $x_1\notin
\q$ because $x_1\notin \n$ and $x_1\in C$.  Since $x_1$ lies in every
maximal ideal of $B$ other than $\m$, necessarily $\q=\m$, as required.
Thus $\beta$ is bijective.

Consider $\alpha$.  It is injective because $C\subseteq B$.  Now,
$\alpha$ is finite; hence, by Nakayama's lemma, it will be surjective if
$\alpha\ox k$ is.  However, $(\beta\alpha)\ox k$ is surjective by
assumption, and $\beta$ is bijective by the preceding paragraph.  Hence
$\alpha$ is surjective, so bijective.  So $C_{\n}\to B_{\m}$ is bijective.
Thus,  to prove (2) as well as (3), we may assume that $R$ is a field.

We may also assume that $B$ is local.  Indeed, since $R$ is a field,
$B$~is equal to the direct product, over its finite set of maximal
ideals $\m$, of its localizations: $B=\prod B_{\m}$.  Suppose one of the
localizations $B_{\m}$ has $t$ generators.  Lift them to elements
$x_1,\dots,x_t$ of $B$ whose images in the other localizations are
equal to 0.  Then form $C$ and $\n$ as usual.  Clearly $C=B_{\m}$ and
$C_{\n}=C$.  Hence, to prove (2), we may replace $B$ by $B_{\m}$, and then we
have to prove that $B$ is generated as an $R$-algebra by $t$~elements.

Consider (3).  Suppose that each $B_{\m}$ has $t$ generators $x_{\m,i}$.
Set $x_i := (x_{\m,i})$ in $B$.  Now, in the polynomial ring in one
variable $R[\lambda]$, let $f_{\m}(\lambda)$ be a polynomial of minimal
degree such that $f_{\m}(x_{\m,1})=0$.  Suppose that $R$ is infinite.  Then
we may assume that the $f_{\m}(\lambda)$ are relatively prime; indeed,
replacing $x_{\m,1}$ by $x_{\m,1}+a_{\m}$ where the $a_{\m}$ are suitable
elements of $R$, we may ensure that no two $f_{\m}(\lambda)$ share a root
in some algebraic closure of $R$.  Set $f:=\prod f_{\m}$.  Then
$R[\lambda]/(f)$ is isomorphic to the subalgebra $B_1$ of $B$ generated
by $x_1$.  Hence, by the Chinese remainder theorem, $B_1$ contains the
idempotents of the decomposition $B=\prod B_{\m}$.  Therefore, $B$ is
generated by the $x_i$.  Thus, to prove (3) as well as (2), we may
assume that $B$ is local.

We may also assume $t=1$.  Indeed, suppose
$\Omega^1_{B/R}/\m\Omega^1_{B/R}$ has dimension at least 2, and let
$x_1$ be an element of $B$ such that the residue class of $dx_1$ is
nonzero.  Let $B_1$ be the subalgebra of $B$ generated by $x_1$.  Then
$B$ is a finitely generated $B_1$-module; so, by the Cohen--Seidenberg
theorem, $B_1$ is local as $B$ is.  Moreover, $\Omega^1_{B/B_1}$ is
equal to the quotient of $\Omega^1_{B/R}$ by the submodule generated by
$dx_1$.  Hence, we may assume by induction on $t$ that there exist
$t-1$ elements $x_2,\dots,x_t$ of $B$ that generate it as a
$B_1$-algebra.  Then $x_1,\dots,x_t$ generate $B$ as an $R$-algebra.

Suppose that the natural map $R\to B/\m$ is an isomorphism.  Then there
exists an element $x$ in $\m$ that generates $B$ as an $R$-algebra.
Indeed, $\m/\m^2$ is equal to $\Omega^1_{B/R}/\m\Omega^1_{B/R}$; see, for
example,  \cite{Kunz, Cor.~6.5(a), p.~96}.   So
$\dim(\m/\m^2)\le1$.  Let $x$ be an element of $\m$ whose residue class
generates $\m/\m^2$ over $R$; possibly, $x=0$.  Then $B=R[x]$.  (The
argument is standard.  Let $y$ be an element of $B$.  For all $n\ge0$,
there is a polynomial $y_n:=\sum_{i=0}^{n}a_ix^i$ with $a_i\in R$ and
$y-y_n$ in $\m^{n+1}$; indeed, take $a_0$ to be the image of $y$ in
$R=B/\m$, and given $y_n$, take $a_{n+1}$ so that $y-y_n$ is equal to
$a_{n+1}x^{n+1}$.  However, $\m^n=0$ for $n\gg0$, so $y=y_n$ for
$n\gg0$.)

Suppose that $R$ is infinite.  Let $R'$ be an algebraic closure of $R$,
and set $B':=B\otimes_RR'$.  Then there exists an element $x'$ of $B'$
that generates it as an $R'$-algebra by the preceding paragraphs,
because $B'$ is a finite $R'$-algebra whose residue class fields are
all equal to $R'$ and because the formation of $\Omega^1$ commutes with
base change.  To descend the existence of a generator, let
$y_1,\dots,y_n$ be a vector space basis of $B$ over $R$, and let
$\tau_1,\dots,\tau_n$ be indeterminates.  In the polynomial ring
$B[\tau_1,\dots,\tau_n]$, set $u:=\sum \tau_iy_i$ and expand the powers $u^j$
for $j=0,\dots,n-1$.  Form the matrix $\Phi(\tau_1,\dots,\tau_n)$ such that
	$$(1,u,\dots,u^{n-1})^{\rm tr}
	 = \Phi(\tau_1,\dots,\tau_n)(y_1,\dots,y_n)^{\rm tr}$$
 where `tr' denotes transpose.  Say $x'=\sum \tau'_iy_i$ where $\tau'_i\in
R'$.  Then the $\tau'_i$ do not satisfy the equation
$\det\Phi(\tau_1,\dots,\tau_n)=0$ because $x'$ generates $B'$ as an
$R'$-algebra.  Hence, since $R$ is infinite, there exist elements
$\tau''_1,\dots,\tau''_n$ in $R$ that do not satisfy this equation.
Therefore, $x:= \sum \tau''_iy_i$ generates $B$ as an $R$-algebra.

Finally, suppose that $R$ is finite, say of characteristic $p$.  Take
$q:=p^e$ so large that $\m^q=0$ and $z^q=z$ for all $z\in B/\m$.
Consider $L:=B^q$.  Obviously, $L$ is an $R$-subalgebra.  In fact, $L$
is a field: if $x\in B$ and $x^q\ne0$, then $x\notin \m$ as $\m^q=0$;
hence, there is a $y\in B$ such that $xy=1$; so $x^qy^q=1$.  Now, since
$z^q=z$ for all $z\in B/\m$, the natural map $L\to B/\m$ is
surjective; so it is an isomorphism because $L$ is a field.  Since
$\Omega ^1_{B/L}$ is a quotient of $\Omega ^1_{B/R}$, by the paragraph
before the last, there is an $x\in \m$ such that $B=L[x]$.  Since $L$ is
finite, its multiplicative group is generated by an element $y$.  Set
$z:=x+y$, and $B':=R[z]$.  Then $z^q=x^q+y^q$, so $z^q=y$.  Thus $y\in
B'$.  Hence $L\subseteq B'$.  Moreover, $z-z^q=x+y-y=x$.  Thus $x\in
B'$.  So $L[x]$ lies in $B'$.  Therefore $B=B'$, and the proof is now
complete.

\dfn9
 Let $f\:X\to Y$ be a map of locally Noetherian schemes.  Following
\cite{SGAVI, 1.1, p.~466} call $f$ a {\it local complete intersection\/}
if, locally on $X$, there is a factorization $f=\pi i$ where
$i\:X\hookrightarrow P$ is a regular embedding and $\pi\:P\to Y$ is
smooth.  Following \cite{HRD, p.~144}, call $f$ {\it Gorenstein\/} if
it has finite flat dimension and if in the derived category $f^!\O_Y$ is
isomorphic to a (shifted) invertible sheaf.

\prp10
 Let $f\:X\to Y$ be a finite map of locally Noetherian schemes, $s$ an
integer.  If $f$ is a local complete intersection and is locally of
codimension $s$, then $f$ is Gorenstein and locally of flat dimension
$s$.  Moreover, the converse holds if also $s=1$ and $f$ is curvilinear.
 \pf
 Suppose $f$ is a local complete intersection and is locally of
codimension $s$.  Then clearly $f$ is Gorenstein; see \cite{HRD,
Cor.~7.3, p.~180, } and \cite{HRD, 3, p.~190}.  Now, for every $x$ in
$X$, clearly
	$$\depth\O_{X,x}=\depth\O_{Y,fx}-s.$$
 Hence, $f$ is locally of flat dimension $s$ by the
Auslander--Buchsbaum formula, which applies after completion because
$f$ is finite.  Thus the direct assertion holds.  The converse follows
from \Cs7)(1) and the following lemma.

\lem11
 Let $R\to B$ be a quasi-finite local homomorphism of Noetherian local
rings, and $t$ an integer.  Assume \smallskip
 \item i that $B$ is of the form $S/I$
where $S$ is a localization at a prime ideal of the polynomial ring in
one variable $R[u]$,
 \item ii that $\Ext_S^i(B,S)$ vanishes for $i\ne t$
and is isomorphic to $B$ for $i=t$, and
 \item iii that $B$ has flat dimension 1 over $R$.
 \smallskip Then $t=2$ and $I$ is generated by a regular sequence of
length $2$; moreover, $\dim B=\dim R-1$.
 \pf
 Assumptions (i) and (iii) imply that $I$ is flat over $R$.  Say $I$ is
generated by $n$ elements with $n$ minimal, and form the exact sequence,
	$$0\To F\To S^n\TO\sigma I\To 0.\dno10.1$$
 Then $F$ is flat over $R$.  Moreover, if $k$ denotes the residue field
of $R$, then $F\ox k$ is free over $S\ox k$, because $S\ox k$ is a
Principal Ideal Domain.  Hence $F$ is a free $S$-module.  So, clearly,
$F=S^{n-1}.$

Sequence \Cs10.1) therefore yields this free resolution of $B$ over $S$:
	$$0\To S^{n-1}\To S^n\To S\To B\To 0.\dno10.2$$
 Hence $t\le2$.  Suppose $t<2$.  Then, because of (ii), dualizing
\Cs10.2) yields a surjection $v\:S^n\onto S^{n-1}$.  However, $v\ox
k=0$ because $n$ is minimal.  Hence $n=1$.  Therefore, the map
$\sigma\:S\to I$ is an isomorphism.  Consequently, there is a short
exact sequence,
	$$0\To \Tor^S_1(B,k)\To S\ox k\To S\ox k\To B\ox k\To 0.$$
 Since $B/R$ is quasi-finite, $\dim_kB\ox k$ is finite.  Hence the map
in the middle is nonzero.  Therefore, it is injective because $S\ox k$
is a domain.  So $\Tor^S_1(B,k)=0$.  Hence $B$ is flat over $R$ by the
local criterion.  Thus  (iii) is contradicted, and so $t=2$.

Because $t=2$ and because of  (ii), dualizing \Cs10.2) yields
the following exact sequence:
	$$0\To S\To S^n\To S^{n-1}\To B\To0.$$
 Because $S$ is local, this sequence may be reduced
to the sequence,
	$$0\To S\To S^2\To S\To B\To0,$$
 by splitting off copies of the trivial exact sequence $0\to S\to
S\to0$.  Therefore, $I$ is generated by two elements, and because
$t=2$, they form a regular sequence.  Finally, by hypothesis, $S$ is
the localization of $R[u]$ at a prime ideal, and this ideal must be
maximal because $B$ is quasi-finite over $R$; hence, $\dim S=\dim R+1$.
Therefore, $\dim B=\dim R-1$ because $I$ is generated by a regular
sequence of length 2.

\newsect	The multiple-point schemes

\dfn1
 A map $f\:X\to Y$ of locally Noetherian schemes will be said to be
{\it birational onto its image\/} if there is an open subset $U$ of $Y$
such that (i) its preimage $f^{-1}U$ is dense in $X$ and (ii) the
restriction $f^{-1}U\to U$ is an embedding.

\prp2
 Let $f\:X\to Y$ be a finite map of locally Noetherian schemes.
 \part1
 The map $f$ is birational onto its image if and only if the source
double-point scheme $M_2$ is nowhere topologically dense in $X$.  These
two equivalent conditions imply that $N_2$ is nowhere topologically
dense in $N_1$, and all three conditions are equivalent if $f$ is
locally of codimension $1$.
 \part2
 The scheme-theoretic image of $f$ is a closed subscheme of the
scheme of target points $N_1$.  The two schemes have the same support,
and they are equal off the scheme of target double-points $N_2$.  If
they are equal everywhere and if $f$ is locally of flat dimension $1$,
then $f$ is birational onto its image.
 \part3
 The map $f$ induces a finite, surjective map $M_r\to N_r$ for $r\ge1$.
 \pf
 Let $U$ be the largest open subset of $Y$ such that $f^{-1}U\to U$
is an embedding.  Since $f$ is finite, $U$ consists of all $y\in Y$ at
which the comorphism $\O_Y\to f_*\O_X$ is surjective.  So, by
Nakayama's lemma, $U$ consists of the $y$ at which the vector space
$(f_*\O_X)(y)$ has dimension at most~1.  Hence, $U=Y-N_2$.  Therefore,
since $M_2=f^{-1}N_2$, the first assertion of (2) holds.  Obviously, if
$M_2$ is nowhere topologically dense in $X$, then $N_2$ is nowhere
topologically dense in $N_1$.  Finally, if $f$ is locally of
codimension 1, then every component of $X$ must map onto a component of
$N_1$; hence, if also $N_2$ is nowhere topologically dense in $N_1$,
then $M_2$ is nowhere topologically dense in $X$.  Thus (1) holds.

 The scheme-theoretic image of $f$ is defined as the smallest closed
subscheme $Z$ of $Y$ through which $f$ factors \cite{EGAI, (6.10.1),
p.~324}.  Because $f$ is quasi-compact and quasi-separated, $Z$ exists
and is associated to the ideal $\ANN_Y(f_*\O_X)$.  Since locally on $Y$
there is an integer $n$ such that
  $$\ANN_Y(f_*\O_X)^n \subseteq \Fit0Y(f_*\O_X) \subseteq\ANN_Y(f_*\O_X),$$
 the image $Z$ is a closed subscheme of $N_1$, and the two schemes have
the same support.  They are equal off $N_2$ because there $f_*\O_X$ is
a cyclic $\O_Y$-module.  If they are equal everywhere, then
$\ANN_Y(f_*\O_X)=\Fit0Y(f_*\O_X)$.  On the other hand, by \cite{B-E,
Thm.~3.1},
	$$\ANN_Y(f_*\O_X)=\Fit0Y(f_*\O_X):\Fit1Y(f_*\O_X).$$
 It follows that the stalk $\Fit1Y(f_*\O_X)_x$ is not contained in any
associated prime of the stalk $\Fit0Y(f_*\O_X)_x$ for any $x\in X$.  So
$N_2$ is nowhere topologically dense in $N_1$.  Therefore, $f$ is
birational onto its image by~(1) and~(2.5)(1).  Thus (2) holds.

By definition, $M_r=f^{-1}N_r$.  Hence $f$ induces a
finite map $M_r\to N_r$.  It is surjective because $N_r\subseteq N_1$
and because $f$ carries $X$ onto $N_1$ by~(1).  Thus (3) holds.

\dfn3
 Following \cite{EGAIV, (5.7.2), p.~103}, a locally Noetherian
scheme $Y$ will be said to satisfy {\it Serre's condition\/} \Ser r, if
for every $y\in Y$,
	$$\hbox{depth}(\O_{Y,y}) \ge \inf(r,\,\dim(\O_{Y,y})).$$

\prp4
 Let $f\:X\to Y$ be a finite map of locally Noetherian schemes.  Assume
that $f$ is locally of flat dimension 1 and is birational onto its
image.  Assume also that $Y$ satisfies \Ser2.  Let $Z$ denote the
scheme-theoretic image of $f$.

Then $Z$ is equal to $N_1$.  Furthermore, $N_2$ is defined by the
adjoint ideal $\ANN_Y(f_*\O_X/\O_{Z})$, and $M_2$ is defined by the
conductor $\C_X$.  Each component of $M_2$ has codimension $1$ and maps
onto a component of $N_2$; each component of $N_2$ has codimension $2$;
and the fundamental cycles of these two schemes are related by the
equation,
	$$f_*[M_2] = 2[N_2].$$
 Finally, $\O_{N_2}$ and $\O_{f_*M_2}$ are perfect $\O_Y$-modules of
grade $2$.
 \pf
 By \Cs2)(2), we have $Z\subseteq N_1$, with equality off $N_2$, and
$N_2$ is nowhere topologically dense in $N_1$ by \Cs2)(1) and (2.5)(1).
Now, $N_1$ has no embedded components because it is a divisor by (2.2)
and because $Y$ satisfies \Ser2.  Hence $N_1=Z$, as asserted.

 The cyclic $\O_Y$-module with ideal $\ANN_Y(f_*\O_X/\O_{Z})$ has a
Hilbert--Burch resolution, and
  $$\ANN_Y(f_*\O_X/\O_{Z})=\Fit1Y(f_*\O_X)=\Fit0Y(f_*\O_X/\O_{Z});\dno4.1$$
 see \cite{KLU, (3.5), p.~208}.  The first equation of \Cs4.1) says
that the adjoint ideal defines $N_2$.  Since $\C_X$ is the ideal on $X$
induced by the adjoint ideal, therefore $\C_X$ defines $M_2$.  Because
of the Hilbert--Burch resolution, \Cs4.1) implies that $\O_{N_2}$ is a
perfect $\O_Y$-module of grade 2.  Hence $\O_{f_*M_2}$ is perfect of
grade 2 too because $Z$ is a divisor.
 Moreover, the determinantal\- length formula
(see \cite{BR, 2.4, 4.3, 4.5} or \cite{BVMZ, 2.10,
p.~154}) gives the following equation, in which $\nu$ is the generic
point of an arbitrary component of $N_2$ and $\ell_{\nu}$ indicates the
length of the stalk at~$\nu$:
  $$\ell_{\nu}(f_*\O_X/\O_{Z})=\ell_{\nu}(\O_Y/\Fit0Y(f_*\O_X/\O_{Z})).$$
 Let $\C_{Z}$ denote the conductor on $Z$, namely, the ideal induced
by the adjoint ideal.  Then, clearly,
  $$f_*(\O_X/\C_X) = (f_*\O_X)/\C_Z
   \hbox{ and }\O_Y/\Fit0Y(f_*\O_X/\O_{Z}) = \O_{Z}/\C_{Z}.$$
 The preceding two displays yield
  $$\ell_{\nu}(f_*(\O_X/\C_X)) = 2\ell_{\nu}(\O_{Z}/\C_{Z}).$$
 Rewritten, the latter equation will become $[N_2] = 2f_*[M_2]$, once we
prove the assertion about the components of $N_2$ and $M_2$.

Let $\eta$ be the generic point of a component of $M_2$.  Then
$\O_{M_2,\eta}$ is an Artin ring, and its residue field is a finite
extension of that of $\O_{Y,f\eta}$.  So $\O_{M_2,\eta}$ is an
$\O_{Y,f\eta}$-module of finite length, and of flat dimension at most~$2$,
hence of projective dimension at most~$2$.  Therefore the
intersection theorem of P.~Roberts \cite{PRob} implies that
$\O_{Y,f\eta}$ has dimension at most~2.  However, every component of
$N_2$ has codimension at least~2.  Consequently, the original component
of $M_2$ has codimension 1 by (2.5)(1).  Finally, every component of
$N_2$ is the image of a component of $M_2$ by (3.2)(3).  The proof is
now complete.

\thm5
 Let $f\:X\to Y$ be a finite map of locally Noetherian schemes, and $r$
an integer, $r\ge1$.  Assume that $f$ is locally of flat dimension $1$
and curvilinear.  Then each component of $N_r$ has codimension at
most~$r$ \(that is, the local ring at the generic point has dimension
at most~$r$).  Assume further that each component of $N_r$ has
codimension~$r$ and that $Y$ satisfies Serre's condition \Ser r.  Then
$\O_{N_r}$ and $\O_{M_r}$ are~perfect $\O_Y$-modules of grade $r$, each
component of $M_r$ has codimension $r-1$ and maps onto a component of
$N_r$, and the fundamental cycles of these two schemes are related by
the equation
	$$f_*[M_r] = r[N_r].$$

 \pf
 The assertion results from (2.2), (3.2)(3) and the next lemma.

 \lem6
 Let $R$ be a local Noetherian ring, and $B$ an $R$-algebra that is
finitely generated as a module.  Assume that, for every maximal ideal
$\m$ of $B$,
	$$\dim_{B/\m}(\Omega^1_{B/R}/\m\Omega^1_{B/R}) \le 1.$$
 Assume that the $R$-module $B$ is presented by a square matrix whose
determinant is regular.  Fix $r \ge1$ and let $F:=\fit R{r-1}(B)$ denote
the Fitting ideal, and $\Ht F$ its height.  If $F\not=R$, then $\Ht
F\le r$.  If $\Ht F = r$ and if $R$ satisfies \Ser r, then both $R/F$
and $B/FB$ are perfect $R$-modules of grade $r$; moreover, then, any
minimal prime of $FB$ has height $r-1$, and its preimage in $R$ is a
minimal prime of $F$ of height $r$.  Finally, if $\dim R = r$ too, then
the following length relation obtains:
	$$\ell_R(B/FB) = r\,\ell_R(R/F).$$
 \pf
 By using a standard device, we may assume that $R$ has an infinite
residue class field: replace $R$ and $B$ by $R'$ and $B\ox_R R'$, where
$R'$ is the flat local $R$-algebra obtained by forming the polynomial
ring in one variable over $R$ and localizing it at the extension of the
maximal ideal of $R$.  Then (2.8)(3) implies that there is an $x$ in
$B$ such that $B=\sum_{i=0}^{n-1}Rx^i$ where $n\ge r$.  By Lemma~(2.3),
there is an $n$ by $n$ matrix $\psi$ with entries in $R$ such that the
corresponding short sequence is exact:
	$$\CD0 @>>> R^n @>\psi>> R^n @>>> B @>>> 0\endCD$$
 where the natural basis element $e_i$ of $R^n$ is mapped to the
generator $x^{i-1}$ of $B$.

Let $M$ be the $R$-module $B/\sum_{i=0}^{r-2}Rx^i$, and let $\phi$ be
the  $n-r+1$ by~$n$  matrix consisting of the last $n-r+1$ rows of
$\psi$.  Then there is a commutative diagram with exact rows and
surjective columns,
 $$\CD
 R^n @>\psi>> R^n=\oplus_{i=1}^nRe_i   @>>> B @>>> 0	\\
     @|		@VVV			@VVV		\\
 R^n @>\phi>> R^{n-r+1}=\oplus_{i=r}^nRe_i @>>> M @>>> 0 \rlap{\quad .}
 \endCD$$
 Let $I_{n-r+1}(\phi)$ and $I_{n-r+1}(\psi)$ denote the ideals of
$n-r+1$ by $n-r+1$ minors. So
  $$F = I_{n-r+1}(\psi) \hbox{  and } \fit R0(M) = I_{n-r+1}(\phi).$$
 Now, Gruson and Peskine \cite{GP, Lem.~1.3, p.~4}
proved (using the multiplicative structure of $B$) that $F =
\fit R0(M)$.  So $F$ is generated by the maximal minors of $\phi$.
So, by the classical height result \cite{BVSpgr, (2.1), p.~10},
	$$\Ht F = \Ht I_{n-r+1}(\phi) \le n-(n-r+1)+1.$$

Suppose $\Ht F =r$ and $R$ satisfies \Ser r.  Then $R/F$ and $M$ are
perfect $R$-modules of grade~$r$ by Eagon's theorem \cite{BVSpgr,
(2.16)(c), p.~18}.

Let `$\mkern1mu\?{\phantom{s}}\mkern1mu$' denote reduction modulo
$I_{n-r+1}(\phi)$.  We will construct a commutative diagram with exact
rows
 $$\def\bwg{\bigwedge^{n-r+1}}	\CD
\R{n-r+1}\ox_R\bwg\Rn@>\?\delta>>\Rn@>\?\psi>>\Rn@>>>B\ox_R\?R@>>>0\\
	@|		   @|			@VVV \\
\R{n-r+1}\ox_R\bwg\Rn@>\?\delta>>\Rn@>\?\phi>>\R{n-r+1}@>>> M @>>>0
 \endCD$$
 in which the right vertical map is surjective.  Once constructed, this
diagram yields an exact sequence,
	$$0\To\R{r-1}\To B\ox_R\?R \To M \To 0.\dno6.1$$
 Now, $B\otimes_R\bar R$ is equal to $B/FB$; hence, $B/FB$ is perfect
over $R$ of grade~$r$ because $\?R$ and $M$ are.

Let $\p$ be a minimal prime of the $B$-ideal $FB$, and let $\q$ be its
preimage in $R$.  Then $(B/FB)_{\p}$ is an Artin ring, and its residue
field is a finite extension of that of $R_{\q}$.  So $(B/FB)_{\p}$ is an
$R_{\q}$-module of finite length, and of flat dimension at most $r$, hence
of projective dimension at most~$r$.  Therefore the intersection
theorem of P.~Roberts \cite{PRob} implies that $R_{\q}$ has dimension at
most $r$.  However, $\Ht F =r$.  Consequently, $\q$ is a minimal prime of
$F$ of height $r$.  So, (2.5)(1) implies that $B_{\p}$ has dimension
$r-1$, as asserted.

If $\dim R =r$ too, then the determinantal length formula \cite{BR,
2.4, 4.3, 4.5} and \cite{BVMZ, 2.10, p.~154} yields
 $$\ell_R(M) = \ell_R(\Cok \phi) = \ell_R(R/I_{n-r+1}(\phi)) = \ell_R(R/F).$$
 Since $R/F=\?R$ and $B\otimes_R\bar R=B/FB$, then \Cs6.1) yields the
asserted length relation.

To define $\?\delta$, use the bases $\{e_r,\,\dots,\,e_n\}$ and
$\{e_1,\,\dots,\,e_n\}$ of $R^{n-r+1}$ and $R^n$.  For $r\le i\le n$,
for $1\le j\le n-r+1$, and for $1\le k_1<\cdots<k_{n-r+1}\le n$, denote
by $d_i^{{\bf k}_j}$ the minor (of $\phi$) that is obtained from $\psi$
by deleting rows 1, \dots, $r-1$, $i$ and taking columns $k_1$, \dots,
$k_{j-1}$, $k_{j+1}$, \dots, $k_{n-r+1}$.  Now, define
	$$\eqalign{&\CD\delta\:R^{n-r+1}\ox_R\bigwedge^{n-r+1}R^n
	 @>>> R^n\qquad \hbox{by}\endCD\cr
	&\delta(e_i\ox e_{k_1}\wedge\cdots\wedge e_{k_{n-r+1}})
	 := \sum\nolimits_{j=1}^{n-r+1}(-1)^jd_i^{{\bf k}_j}\cdot e_{k_j}. \cr}$$
 It is easy to see that the image of the composite map,
 $$\CD
 R^{n-r+1}\ox_R\bigwedge^{n-r+1}R^n @>\delta>> R^n @>\phi>>R^{n-r+1}
 \endCD,$$
 is exactly $I_{n-r+1}(\phi)\cdot R^{n-r+1}$.

On the other hand, since $I_{n-r+1}(\phi)$ has generic grade in $R$, it
follows that an $R$-resolution of $M=\Cok \phi$   is given by the
Buchsbaum--Rim complex   \cite{BR, 2.4, p.~207},
 $$\CD
 \cdots@>>>\bigwedge^{n-r+2}R^n@>\beta>>R^n@>\phi>>R^{n-r+1}@>>>M@>>>0
 \endCD$$
 where $\?\beta=0$.  Since $\phi$ maps $\Im\delta$ onto
$I_{n-r+1}(\phi)\cdot R^{n-r+1}$, the preimage
$\phi^{-1}(I_{n-r+1}(\phi)\cdot R^{n-r+1})$  is exactly
$\Im\delta+\Im\beta$.  Therefore, the sequence
  $$\CD\matrix{\bigwedge^{n-r+2}\Rn\cr \bigoplus\cr
		\R{n-r+1}\otimes\bigwedge^{n-r+1}\Rn\cr}
  @>\?\beta\oplus\?\delta>>\Rn @>\?\phi>>\R{n-r+1}@>>>M\ox_R\?R@>>>0
  \endCD$$
 is exact.  Thus the bottom row of the diagram is exact, because
$\?\beta=0$ and  $M\ox_R\?R=M$.

It is also easy to see that the image of the composite map
 $$\CD R^{n-r+1}\ox_R\bigwedge^{n-r+1}R^n @>\delta>> R^n @>\psi>> R^n\endCD$$
 is contained in $I_{n-r+1}(\psi)\cdot R^n$, which is equal to
$I_{n-r+1}(\phi)\cdot R^n$.  Therefore, the top row of the diagram
is a complex.  This complex is exact, because every relation on
the columns of $\?\psi$ is a relation on the columns of $\?\phi$ and
hence contained in the image of $\?\delta$.  The proof is now complete.

\cor7
 Let $f\:X\to Y$ be a finite map of locally Noetherian schemes, and $r$
an integer, $r\ge1$.  Assume that $f$ is locally of flat dimension 1,
and, if $r\ge3$, then assume that $f$ is curvilinear.  Assume also that
$Y$ satisfies \Ser {r+1}, that each component of $N_r$ has codimension~$r$,
and that each component of $N_{r+1}$ has codimension $r+1$.  Then
$N_r$ is the scheme-theoretic image of $M_r$.
 \pf
 Since $M_r=f^{-1}N_r$ and $f$ is finite, $f_*\O_{M_r}$ is equal to the
restriction of $f_*\O_X$ to $N_r$.  Hence, $f_*\O_{M_r}$ is locally
free of rank $r$ on $N_r-N_{r+1}$ by standard linear algebra.  So the
comorphism $\gamma \:\O_{N_r}\to f_*\O_{M_r}$ is injective off
$N_{r+1}$.  Now, $N_r$ is perfect by (2.2) if $r=1$; by \Cs2)(2),
(2.5)(1), and \Cs4) if $r=2$; and by \Cs5) if $r\ge3$.  Hence
$\O_{N_r}$ has no embedded points because $Y$ satisfies \Ser {r+1}.
Therefore, $\gamma$ is injective everywhere because each component of
$N_{r+1}$ has codimension $r+1$.  The proof is now complete.

\dfn8 Let $f\:X\to Y$ be a finite map of locally Noetherian schemes.
Following \cite{KActa, 4.1, pp.~36--37}, \cite{KSitges, (2.10),
\pp.112--113}, and \cite{KLU, (3.1)}, define the {\it iteration scheme}
$X_2$ and the {\it iteration map} $f_1\:X_2\to X$ of $f$ as follows:
	$$X_2:=\IP(\Idiag)\quad \hbox{and}\quad
	  f_1\:X_2\mathop{\To}^p X\times_YX\mathop{\To}^{p_2}X,$$
 where $\Delta$ is the diagonal, $\Idiag$ is its ideal, $p$ is the
structure map, and $p_2$ is the second projection.  (Thus, $X_2$ is the
residual scheme of $\Delta$.)

\lem9
 Let $f\:X\to Y$ be a finite map of locally Noetherian schemes, and
assume that $f$ is curvilinear.  Then, for any $r\ge2$,
	$$M_r(f) = N_{r-1}(f_1).$$
 \pf
 The structure map $p:X_2\to X\times_YX$ is a closed embedding if and
only if $f$ is curvilinear.  If so, then $p_*\O_{X_2}$ is locally
isomorphic to $\Idiag$, and therefore, for any $r\ge0$,
	$$\Fit rX(f_{1*}\O_{X_2})=\Fit rX(p_{2*}{\Idiag}).\dno9.1$$
 These statements are not hard to prove; see \cite{KLU, (3.2), (3.4)}.
In that  reference, \Cs9.1) is stated only for $r=0$, but
the proof works without change for any $r$.

Since $f$ is an affine map, the operator $f_*$ is exact and
commutes with base change.  Hence, applying  $p_{2*}$ to the natural exact
sequence,
	$$0 \To \Idiag \To \O_{X\times X} \To \O_\Delta \To 0,$$
yields an exact sequence,
	$$0\To p_{2*}{\Idiag} \To f^*f_*\O_X \To \O_X \To 0.$$
 Hence, by standard properties of Fitting ideals,
	$$\Fit rX(p_{2*}{\Idiag}) = \Fit{r+1}X(f^*f_*\O_X)
	 = \Fit{r+1}Y(f_*\O_X)\,\O_X.$$
 Therefore, \Cs9.1) yields the assertion.

\lem10
 Let $f\:X\to Y$ be a finite and curvilinear map of locally Noetherian
schemes.
 \part1 Then $f_1\:X_2\to X$ is finite and curvilinear.
 \part2 Assume that $X$ has no embedded components.  Assume either \(i)
that each component of $N_2$ has codimension at least $2$ in $Y$, or \(ii)
that each component of $M_2$ has codimension at least $1$ in $X$, or
\(iii) that $f$ is birational onto its image.  Finally, assume that $f$
is a local complete intersection and is locally of codimension $1$.  Then
$f_1$ is also a local complete intersection and locally of codimension
$1$.
 \pf
 Assertion (1) holds because (a) the map $p:X_2\to X\times_YX$ is a
closed embedding; (b) the projection $p_2\:X\times_YX \to X$ is finite;
and (c) $\Omega ^1_{p_2}=p_2^*\Omega ^1_f$.

Consider (2).  Conditions (i) and (ii) are equivalent because $f$
is locally of codimension $1$.  Conditions (ii) and (iii) are
equivalent by \Cs2)(1); in particular, (ii) obtains.  By (1),
$f_1\:X_2\to X$ is finite, and by \Cs2)(2), it factors through
$N_1(f_1)$.  By \Cs9), $N_1(f_1)=M_2$.
So, for any $x\in X_2$,
	$$\dim\O_{X_2,x}\le\dim\O_{N_1(f_1),f_1x}
	=\dim\O_{M_2,f_1x}\le\dim\O_{X,f_1x}-1.$$
 Since $X$ has no embedded components, and since $f$ is a local
complete intersection and locally of codimension 1, it follows that
$f_1$ is also.  Indeed, it is not hard to show, see \cite{KActa, 4.3,
p.\thinspace 39}, that $X_2$ is, locally at any point $x$, cut out of some
smooth
$X$-scheme $P$, say of relative dimension~$p$, by $p+1$ elements.
Since $f_1$ is finite and locally of codimension~1, a subset of $p$ of
the elements must restrict to a system of parameters in the fiber of
$P$ through $x$.  Since the fiber is smooth, this system is a regular
sequence.  Hence, by the local criterion of flatness, the $p$ elements
themselves form a regular system and they cut out of $P$ a flat
$X$-scheme~$Q$.  Since $X$ has no embedded components, neither does $Q$
because $Q$ is $X$-flat.  Since the remaining element cuts $X_2$ out of
$Q$, it is regular on $Q$.  Thus the $p+1$ elements form a regular
sequence on $P$.  Thus (2) holds.

\thm11
 Let $f\:X\to Y$ be a finite map of locally Noetherian schemes, and $r$
an integer, $r\ge2 $.  Assume that $f$ is a local complete
intersection, is locally of codimension $1$, and is curvilinear.
Assume also that either \(i) each component of $N_2$ has codimension at
least $2$ in~$Y$, or \(ii) each component of $M_2$ has codimension at
least $1$ in $X$, or \(iii) $f$ is birational onto its image.

If $X$ has no embedded components, then each component of $M_r$ has
codimension at most $r-1$ in $X$.  Furthermore, if each component of
$M_r$ has codimension $r-1$ and if $Y$ satisfies Serre's condition \Ser
r, then $M_r$ is a perfect subscheme of $X$.
 \pf
 If $Y$ satisfies \Ser r, then $X$ satisfies \Ser{r-1} because $f$ is
locally of codimension 1, and because, for every $x$ in $X$, clearly
	$$\depth\O_{X,x}=\depth\O_{Y,fx}-1$$
 since $f$ is also a local complete intersection.  Hence, in any event,
$X$ has no embedded components.  Therefore, $f_1\:X_2\to X$ is a local
complete intersection and locally of codimension 1 by \Cs10)(2).
Moreover, $f_1$ is finite and curvilinear by \Cs10)(1).  Hence $f_1$
is locally of flat dimension~1 by (2.10).  Therefore, \Cs9) and \Cs5)
yield the assertions.

\prp12
 Let $f\:X\to Y$ be a finite map of locally Noetherian schemes.
Assume that $f$ is curvilinear and that $Y$ satisfies \Ser2.  Then the
following conditions are equivalent:\smallbreak
 \item i $f$ is a local complete intersection, is locally of
codimension~$1$, and is birational onto its image;
 \item ii $f$ and $f_1$ are both locally of flat dimension $1$;
 \item iii $f$ is locally of flat dimension $1$ and $M_2$ is a divisor;
 \item iv $f$ is locally of flat dimension $1$, is birational onto its
image, and is Gorenstein.
 \smallskip
 \noindent Moreover, in \(i) or \(iv) or both, the condition
that $f$ is birational onto its image may be replaced either by the
condition that $M_2$ is nowhere topologically dense in $X$ or by the
condition that $N_2$ is nowhere topologically dense in $N_1$.
 \pf
 In the course of proving \Cs11), it was shown that (i) implies (ii).
Assume (ii).  Then $N_1(f_1)$ is a divisor by (2.2), and $N_1(f_1)=M_2$
by \Cs9).  Thus (iii) holds.  Next, assume (iii).  Then $f$ is
birational onto its image by \Cs2)(2), and it is locally of codimension
1 by (2.5)(1).  Moreover, the ideal of $M_2$ is equal to the conductor
$\C_X$ by \Cs4); so $\C_X$ is invertible.  Hence $f$ is Gorenstein by
\cite{KLU, (2.3)}.  Thus (iv) holds.  Now, (iv) implies (i) by (2.10).
Finally, the last assertion holds by \Cs2)(2).

\newsect	The Hilbert scheme

\MinNo=5   

\prp1
 Let $f\:X\to Y$ be a finite map of locally Noetherian
schemes, and assume that $f$ is curvilinear.  Let $r\ge2$.  Then the
universal subscheme $\Univ^r_f$ of $\Hilb^r_f\times_Y X$ is equal to
the Hilbert scheme  $\Hilb^{r-1}_{f_1}$ of the iteration map
$f_1\:X_2\to X$ defined in (3.8),
	$$\Univ^r_f = \Hilb^{r-1}_{f_1}. $$

\pf For convenience, set $\U^r_f:=\Univ^r_f$ and $\H^r_f:=\Hilb^r_f$.
It will be shown that both $\U^r_f$ and $\H^{r-1}_{f_1}$ have canonical
closed embeddings in $\H^{r-1}_{f}\times_YX$ and then that the two
subschemes are equal.  First of all, the structure map $p\:X_2\to
X\times_YX$ is a closed embedding because $f$ is curvilinear.  Hence,
there is a canonical embedding of $\H^{r-1}_{f_1}$ in $\H^{r-1}_{p_2}$,
which is equal to $\H^{r-1}_{f}\times_YX$.

Secondly, there is a canonical map $v:V\to\U^r_f$ where $V$ is the
residual scheme of $\U^{r-1}_f$ in $\H^{r-1}_{f}\times_YX$ by
\cite{KSitges, (2.9)(1), p.~111}.  The map $v$ is an isomorphism by
\cite{KSitges, (2.9)(4), p.~111}; indeed, every length-$r$ subscheme
$z$ of every fiber $f^{-1}(y)$ is Gorenstein, because $f^{-1}(y)$ is
isomorphic, locally at each of its points, to a closed subscheme of the
affine line over the field $k(y)$ by (2.7)(2) as $f$ is curvilinear.
Moreover, the structure map $V\to \H^{r-1}_{f}\times_YX$ is a closed
embedding because the ideal of $\U^{r-1}_f$ is locally generated by a
single element.  Indeed, the formation of this ideal commutes with base
change through $\H^{r-1}_{f}$ because $\U^{r-1}_f$ is flat, and on each
fiber of $\H^{r-1}_{f}\times X$, the ideal is generated by a single
element; the latter obtains because the fiber comes via base field
extension from a fiber $f^{-1}(y)$, and as noted above, $f^{-1}(y)$ is
isomorphic to a closed subscheme of the affine line over the field
$k(y)$.

Consider the $r$th iteration scheme $X_r$ and the corresponding
iteration map $f_{r-1}\:X_r \to X_{r-1}$ of $f$.  For $r=2$, they are
simply the iteration scheme $X_2$ and the iteration map $f_1\:X_2\to
X$, and, for $r\ge3$, they are defined recursively as the iteration
scheme and iteration map of $f_{r-2}$; see \cite{KActa, 4.1,
p.~36--37} or \cite{KSitges, (4.4), p.~120}.  Since $f$ is
curvilinear, there is a canonical finite, flat, and surjective map
$u\:X_r\to \U^r_f$ by \cite{KSitges, (5.10)(i), p.~128}.  Then
$u_*\O_{X_r}$ is a locally free $\O_{\U^r_f}$-module, so the comorphism
$\O_{\U^r_f}\to u_*\O_{X_r}$ is injective.  Hence, $\U^r_f$ is equal to
the scheme-theoretic image of $X_r$ in $\H^{r-1}_{f}\times_YX$.

It is clear from the definition of $X_r$ that it is equal to the
$(r-1)$st iteration scheme of $f_1$.  Hence, there is a canonical
finite, flat, and surjective map $X_r\to \H^{r-1}_{f_1}$ by
\cite{KSitges, (5.10)(i), p.~128}.  This map yields a second map from
$X_r$ to $\H^{r-1}_{f}\times_YX$, and its scheme-theoretic image is
equal to~$\H^{r-1}_{f_1}$.  It may be checked using the universal
property of the Hilbert scheme that the two maps from $X_r$ to
$\H^{r-1}_{f}\times_YX$ are equal.  Therefore, $\U^r_f$ and
$\H^{r-1}_{f_1}$ are equal too.

\thm2
 Let $f\:X\to Y$ be a finite map of locally Noetherian schemes, and
$r\ge1$.  Assume that $f\:X\to Y$ is a local complete intersection,
locally of codimension $1$, and curvilinear.  Assume that $Y$ satisfies
Serre's condition \Ser{r+1}.  Finally, assume that each component of
$N_s$ has codimension $s$ for $s=1,\dots,r+1$.  Let $h\:\Hilb^r_f\to Y$
denote the structure map.
 Then $h$ is finite, locally of flat dimension $r$, locally of
codimension $r$, and Gorenstein.  Moreover, $h^{-1}N_{r+1}$ is a
divisor, and
	$$h_*[h^{-1}N_{r+1}] = (r+1)[N_{r+1}].$$
 Similar assertions hold for the structure map $h_1\:\Univ^r_f\to X$
too.
 \pf
 First of all, $h$ has finite fibers because $f$ does.  Hence, $h$
is finite because it is proper.  Now, $f$ is locally of flat dimension
1 by (2.10).  Hence, by (3.5),
	$$f_*[M_{r+1}] = (r+1)[N_{r+1}].\dno2.1$$

The proof proceeds by induction on $r$.  Suppose $r=1$.  Then $h$ is
equal to $f$, and the asserted equation becomes \Cs2.1).  Now, $N_2$ is
nowhere topologically dense in $N_1$; hence, $f$ is locally of flat
dimension 1, locally of codimension 1, and Gorenstein, and $M_2$ is a
divisor by (3.12).  However, $M_2=f^{-1}N_2$ essentially by definition.
Thus the assertions about $h$ hold.  Furthermore, $\Univ^1_f$ is equal
to the diagonal subscheme of $X\times_Y X$.  Hence the assertions about
$h_1$ hold too when $r=1$.

Suppose $r\ge2$.  Consider the map $f_1\:X_2\to X$ and the diagram
 $$\CD
 X @<h_1<< \Univ^r_f @= \Hilb^{r-1}_{f_1} \\
@V fVV     @V uVV \\
 Y @<h<<   \Hilb^r_f
 \endCD$$
 in which $h_1$ and $u$ are the natural maps and the equality is that
of \Cs1).  Since $f$ is a local complete intersection and is locally of
codimension~1 and since $Y$ satisfies \Ser{r+1}, clearly $X$ satisfies
\Ser{r}.  In particular, $X$ has no embedded components.  So $f_1$ is
finite, curvilinear, a local complete intersection, and locally of
codimension~1 by (3.10).  Now, $N_{s}(f_1)=\nobreak M_{s+1}$ for $s\ge1$
by (3.9), and $f$ induces a finite, surjective map $M_{s+1}\to N_{s+1}$
by (3.7); hence, $N_s(f_1)$ is of pure codimension~$s$ for
$s=1,\dots,r$.

The induction hypothesis therefore applies to $f_1$.  Hence, $h_1$ is
locally of flat dimension $r-1$, locally of codimension~$r-1$, and
Gorenstein; moreover, $h_1^{-1}M_{r+1}$ is a divisor, and
	$$h_{1*}[h_1^{-1}M_{r+1}] = r[M_{r+1}].$$
 Since $f$ is locally of flat dimension~1, locally of codimension~1,
and Gorenstein, therefore $fh_1$ is locally of flat dimension~$r$,
locally of co\-dimension~$r$, and Gorenstein.  (With the residue field of
an arbitrary point in the image of $fh_1$ as first argument, the
``change of rings'' spectral sequence for `Tor' shows that $fh_1$ is
locally of flat dimension at least~$r$.)

Since $fh_1$ is locally of flat dimension~$r$ and locally of
co\-dimension~$r$, so is $h$ because $fh_1=hu$ and because $u$ is flat
and finite.  Also because $u$ is finite, the dualizing complexes of
$hu$ and $h$ are related by the formula,
	$$u_*\omega_{hu}=\Hom(u_*\O_{\Univ^r_f},\ \omega_h).$$
 Since $u_*\O_{\Univ^r_f}$ is locally free and since $hu$ is Gorenstein
(being equal to~$fh_1$), it follows that $h$ is Gorenstein.  Finally,
since $h_1^{-1}M_{r+1}$ is a divisor, so is $h^{-1}N_{r+1}$ because
$fh_1=hu$ and because $u$ is flat and finite.  Since $u$ is of degree $r$,
	$$u_*[h_1^{-1}M_{r+1}] = r[h^{-1}N_{r+1}].$$
Since $f_*h_{1*} = h_*u_*$,  therefore
	$$h_*[h^{-1}N_{r+1}] = f_*[M_{r+1}].$$
Consequently, the asserted equation follows from \Cs2.1).  Thus, the
theorem is proved.

\thm3
 Under the conditions of \Cs2), the Hilbert scheme $\Hilb^r_f$ is equal
to the blowup $\Bl(N_r,N_{r+1})$, and the universal subscheme
$\Univ^r_f$ of $\Hilb^r_f\times_Y X$ is equal to the blowup
$\Bl(M_r,M_{r+1})$; that is,
	$$\Hilb^r_f = \Bl(N_r,N_{r+1}) \hbox{ and }
	  \Univ^r_f = \Bl(M_r,M_{r+1}).$$
 \pf
 First of all, the structure map $h\:\Hilb^r_f\to Y$ factors
through a map
	$$\beta\:\Hilb^r_f\To\Bl(N_r,N_{r+1}),$$
 which restricts to an isomorphism off $h^{-1}N_{r+1}$.  Indeed, a map
$g\:G\to Y$ factors through $N_r-N_{r+1}$ if and only if $g^*f_*\O_X$
is locally free of rank~$r$ by \cite{Mum, (*), p.~56}.  Hence, $h$
induces an isomorphism,
	$$(\Hilb^r_f-h^{-1}N_{r+1}) \risom (N_r-N_{r+1}) .$$
 Therefore, the ideal of $h^{-1}N_r$ in $\Hilb^r_f$ vanishes off
$h^{-1}N_{r+1}$.  So the ideal vanishes everywhere because
$h^{-1}N_{r+1}$ is a divisor by \Cs2).  Consequently, $h$ factors
through $N_r$.  Therefore, since $h^{-1}N_{r+1}$ is a divisor, the
universal property of the blowup implies that $h$ factors through a map
$\beta $, as claimed.

For convenience, set $B:=\Bl(N_r,N_{r+1})$, denote the exceptional
divisor by $E$, and set $U:=B-E$.  To construct an inverse $\gamma $ to
$\beta $, it suffices to construct a length-$r$ subscheme $Z$ of
$X\times B/B$ whose restriction over~$U$ is equal to $X\times U$.
Indeed, such a $Z$ defines a map $\gamma \:B\to \Hilb^r_f$ such that
$\beta\gamma $ is equal to the identity off $E$ and $\gamma \beta $ is
equal to the identity off $h^{-1}N_{r+1}$.  Since $E$ and
$h^{-1}N_{r+1}$ are both divisors and since $B$ and $\Hilb^r_f$ are
both separated over $N_r$, each composition is equal to the identity
everywhere.  (Indeed, each is equal to the identity on a closed
subscheme of the source because its target is separated; this subscheme
is equal to the source because it contains an open subscheme that
includes every associated point of the source.)

Let $\iota \:U\to B$ denote the inclusion, let $f_B$ and $f_U$ denote
the base extensions of $f$, and let ${\cal E}$ denote the image of
$f_{B*}\O_{X\times B}$ in $\iota _*(f_{U*}\O_{X\times U})$.  Since
${\cal E}$ is the image of an $f_{B*}\O_{X\times B}$-map, ${\cal E}$ is
an $f_{B*}\O_{X\times B}$-module.  Hence ${\cal E}$ is equal to the
direct image of the structure sheaf of a subscheme $Z$ of $X\times
B/B$.  This $Z$ has the desired properties, because ${\cal E}$ is
locally free of rank $r$, as will now be proved.

The question is local on $B$.  Now, each point of $B$ has a
neighborhood~$V$ on which $f_{B*}\O_{X\times B}$ has a free quotient
${\cal F}$ of rank $r$ by \Cn2)(3) applied to any matrix $\bf X$ presenting
$f_{B*}\O_{X\times B}$ over the local ring $R$ of the point and applied
with any minor generating the $(r+1)$st Fitting ideal as $\D_{\bf i}$.
On $U\cap V$, the canonical surjection from $f_{B*}\O_{X\times B}|V$ to
${\cal F}$ is an isomorphism because the source is locally free of rank
$r$.  Hence there is an induced map $u\:{\cal F}\to{\cal E}|V$, which
is an isomorphism on $U\cap V$.  Since $E$ is a divisor, ${\cal F}$ has no
associated point off $U$.  Hence, $u$ is injective on all of $V$.  On
the other hand, $u$ is surjective because ${\cal E}$ is a quotient of
$f_{B*}\O_{X\times B}$.  Thus ${\cal E}$ is locally free, and the first
assertion is proved.

The second assertion follows from the first applied to $f_1\:X_2\to X$
because of \Cs1) and because $f_1$ satisfies the corresponding
hypotheses; the claim about $f_1$ was established in the proof of \Cs2).

\thm4
 Under the conditions of \Cs2), the structure map $h\:\Hilb^r_f\to N_r$
is finite and birational, its conductor is equal to the ideal ${\cal
J}_r$ of $N_{r+1}$ in $N_r$, and reciprocally, $h_*\O_{\Hilb^r_f}$ is
equal to $\Hom({\cal J}_r,\O_{N_{r}})$.  Moreover, ${\cal J}_r$ is
locally a self-linked ideal of $\O_{N_r}$; in fact, locally there exist
sections $t$ of ${\cal J}_r$ such that ${\cal J}_r\O_{\Hilb^r_f}$ is
equal to $t\O_{\Hilb^r_f}$, and ${\cal J}_r=(t\O_{N_r}):{\cal J}_r$ for
any such $t$.  Furthermore, if $r\ge2$, then similar assertions hold
for the structure map $h_1\:\Univ^r_f\to M_r$.
 \pf
 First of all, the assertions about $h_1$ follow formally from those
about $h$; see the third paragraph of the proof of (4.2).  Now, $h$ is
finite and birational by \Cs2) and \Cs3), and $h^{-1}{\cal J}_r$ is
invertible by \Cs2).  Hence, locally, $h^{-1}{\cal J}_r$ is generated
by a single section of ${\cal J}_r$.  The remaining three assertions
are local on $Y$; so we may assume that $Y$ is the spectrum of a local
ring $R$.

By using the standard device of making a suitable (faithfully) flat
change of base, we may assume that $R$ has an infinite residue class
field; namely, we may, clearly, replace $R$ by the flat local
$R$-algebra obtained by forming the polynomial ring in one variable
over $R$ and localizing it at the extension of the maximal ideal of
$R$.  Then $f_*\O_X$ can be presented by a square matrix ${\bf X}$ that
satisfies the hypotheses of (5.9) below; indeed, the condition on
$\grade I_i({\bf X})$ follows from the hypotheses, and the condition
$I_i({\bf X})=I_i({\bf X}_i)$ follows by the reasoning in the first two
paragraphs of the proof of (3.6).
 Hence (5.9) implies that, in the
local ring $A$ of $N_r$, there are an $A$-regular element $\D$ and an
ideal $I$ containing $\D$ such that $IJ=\D J$ and $J=(\D):I$ where $J$
is the ideal in~$A$ of $N_{r+1}$.  Hence, \Cs5) will yield the
remaining three assertions after we prove that the Hilbert scheme
$\Hilb^r_f$ and the two blowups $\Bl(I)$ and~$\Bl(J)$ are all equal.

The isomorphism $\gamma$ in the proof of \Cs3) clearly factors as
follows:
	$$\CD\gamma\:\Bl(J)@>\eta>>\Bl(I)@>\theta>>\Hilb^r_f\endCD$$
 where $\eta$ is the map given by \Cs5)(5) and $\theta$ is given by a
construction similar to that of $\gamma$, but based on the fact that
$I$ is generated by elements of the form given in (5.9).  The
composition $\eta\gamma^{-1}\theta$ is equal to the identity off the
exceptional divisor of $\Bl(I)$; so it is equal to the identity
everywhere, because $\Bl(I)$ is separated over $N_r$.  Therefore, the
maps $\eta$ and $\theta$ are isomorphisms, and the proof is complete.

\def\iod{{\scriptstyle I\over \scriptstyle \Delta}}
\def\dci{(\Delta):I}
 \lem5
 Let $A$ be a ring, $\Delta $ an $A$-regular element, $I$ an ideal
containing $\Delta $, and $J$ an ideal containing $I$.  Let $K$ be the
total quotient ring of $A$. Set
 $B:=A\left[\iod\right]$ and $C:=\{x\in K | xB \subset A \}$.
 \part1 If $J=\dci$, then $C\subset J$.
 \part2 If $IJ=\Delta J$, then $JB=J$ and $J\subset C$.
 \part3 If $IJ=\Delta J$ and if $JB$ is invertible, then $B=\{x\in
K|xJ\subset J\}$.  If in addition $J=\dci$, then $B=\{x\in K|xJ\subset
A\}$.
 \part4 If $IJ=\Delta J$ and if $J$ is finitely generated, then
$\Spec(B)=\Bl(I)$.
 \part5 If $IJ=\Delta J$, then there is an $A$-map
$\eta\:\Bl(J)\to\Bl(I)$.
 \part 6 If $J=C$ and if $J=tB$ for some $t$, then $t$ is an
$A$-regular element of $A$, and $J=tA:J$.
 \pf
 (1) Let $x\in C$.  Then $x=x\cdot 1$, so $x\in A$.  Moreover,
$x(\iod)\subset A$, so $xI\subset \Delta A$.  Hence $x\in\dci$, but
$\dci=J$.

(2) By hypothesis, $IJ=\Delta J$.  So $J(\iod)=J$.  Hence,
$J(\iod)^n=J$ for any $n\ge1$.  Therefore, $JB=J$.  Consequently,
$J\subset C$.

(3) Let $x\in K$, and suppose $xJ\subset J$.  Then $xJB\subset JB$.
Hence $x\in B$ because $JB$ is invertible.  Conversely, if $x\in B$,
then $xJ\subset J$ by (2).

Let $y\in K$, and suppose $yJ\subset A$.  Then $yJB\subset A$ by (2).
So $yJ\subset C$ by definition of $C$.  If in addition  $J=\dci$, then
$yJ\subset J$  by (1), and so $y\in B$ by the preceding paragraph.

(4) First, consider any local $A$-algebra $D$ such that $ID$ is
invertible.  Say $ID=dD$ and $\Delta=ed$.  Then $IJD=\Delta JD$.  So
$JD=eJD$ because $d$ is regular on $D$.  Hence $e$ is a unit by
Nakayama's lemma because $J$ is finitely generated.  Hence $ID=\Delta
D$.  Therefore, the map  $A\to D$  factors through $B$.

Obviously, $\Spec(B)$ is a principal open subscheme of $\Bl(I)$.  Let
$x\in\Bl(I)$ and set $D:=\O_x$.  By the preceding observation, there is
an $A$-map from $\Spec(D)$ to $\Spec(B)$.  This map agrees with the
canonical map of $\Spec(D)$ into $\Bl(I)$ because $\Bl(I)$ is separated
and the two maps agree off the closed subscheme $V(ID)$, which is a
divisor.  Hence, $x\in\Spec(B)$.

(5) Since $IJ=\D J$ and since $J\O_{\Bl(J)}$ is invertible,
$I\O_{\Bl(J)}$ is generated by $\D$.  Moreover, $\D$ is regular on
$\O_{\Bl(J)}$ because it is regular on the complement of the
exceptional divisor.  Thus $I\O_{\Bl(J)}$ is invertible.  Hence the
asserted map $\eta$ exists.

(6) Since $t\in J$, also $t\in A$.  Since $\Delta=bt$ for some $b\in B$
and since $\Delta$ is $A$-regular, so is $t$.  Now, $J^2=Jt$ because
$J=tB$; hence, $J\subseteq tA:J$.  Finally, suppose $x\in tA:J$.
Then $xtB\subset At$.  Hence $x\in C$, but $C=J$.  Thus
$J\supseteq tA:J$, and the proof is complete.

\lem\Mn1
 Let $R$ be a ring, and ${\bf X}$ an $m$ by $n$ matrix.  Fix $p\ge1$,
and set $A:=R/I_{p+1}({\bf X})$ and $J:=I_{p}({\bf X})A$ where
$I_{q}({\bf X})$ denotes the ideal of $q$ by $q$ minors.  Denote the
image in $J$ of the minor of ${\bf X}$ formed using rows
$i_1,\dots,i_p$ and columns $k_1,\dots,k_p$ by $d_{\bf i}^{\bf k}$.
 \part1 Let ${\bf R}_i$ denote the $i$th row of ${\bf X}$.  Then, for any
${\bf i}$ and ${\bf k}$,
 $$d_{\bf i}^{\bf k}{\bf R}_i
  =\sum\nolimits_{j=1}^p(-1)^{j+p}d_{{\bf i}i_j}^{\bf k}{\bf R}_{i_j}$$
 where ${\bf i}i_j$ is the sequence $i_1,\dots,i_p,i$ without its $j$th
element.
 \part2 (Sylvester's relation)
 Then $d_{\bf i}^{\bf k}d_{\bf j}^{\bf l} =d_{\bf j}^{\bf k}d_{\bf
i}^{\bf l}$ for any ${\bf i}$, ${\bf j}$, ${\bf k}$, and ${\bf l}$.
 \pf
 We may assume that ${\bf X}$ is a matrix of indeterminates and that
$R$ is obtained by adjoining them to the integers.  Then $A$ is a
domain.

To prove (1), form a $p+1$ by $n$ matrix ${\bf Y}$ using rows ${\bf
R}_{i_1}$, \dots, ${\bf R}_{i_p}$, ${\bf R}_{i}$.  For each $k$, form a
$p+1$ by $p+1$ matrix ${\bf Y}^{(k)}$ by taking out of ${\bf Y}$
columns $k_1,\dots,k_p$ and column $k$.  Finally, expand the
determinant of ${\bf Y}^{(k)}$ along the last column to get the
asserted equation.

To prove (2), denote the $p$ by $n$ submatrix of ${\bf X}$ consisting
of rows $i_1,\dots,i_p$ by ${\bf X}_{\bf i}$.  Set $d:=d_{\bf j}^{\bf
l}$.  Then, (1) implies that there is a $p$ by $p$ matrix $\bf M$ such that
$d{\bf X}_{\bf i}=\bf M{\bf X}_{\bf j}$; here $\bf M$ depends on ${\bf i}$,
${\bf j}$ and ${\bf l}$, but not on ${\bf k}$.  Hence,
	$$d^pd_{\bf i}^{\bf k}d_{\bf j}^{\bf l}
	=|\bf M|d_{\bf j}^{\bf k}d_{\bf j}^{\bf l}
	=d^pd_{\bf j}^{\bf k}d_{\bf i}^{\bf l}.$$
 Since  $d\ne0$ and $A$ is a domain, the assertion follows.

\lem\Mn2
 Preserve the conditions of \Cn1).  Let ${\bf k}$ range over all
sequences $1\le k_1<\dots<k_p\le n$.  Given elements $a^{\bf k}$ of
$A$, set
	$$\D_{\bf i}:=\sum\nolimits_{\bf k} a^{\bf k}d_{\bf i}^{\bf k}.$$
 Let $I$ be the ideal generated by the various $\D_{\bf i}$.
 \part1
 Let ${\bf X}^{\bf l}$ be the submatrix of ${\bf X}$ consisting of
columns $l_1,\dots,l_p$.  Then $d_{\bf j}^{\bf l}I=\D_{\bf j}I_p({\bf
X}^{\bf l})$.
 \part2
 Let ${\bf X}_{\bf j}$ be the submatrix of ${\bf X}$ consisting of rows
$j_1,\dots,j_p$.  Then $\D_{\bf i}I_p({\bf X}_{\bf j})\subseteq\D_{\bf j}J$.
Furthermore, if $J=I_p({\bf X}_{\bf j})A$, then $IJ=\D_{\bf j}J$.
 \part3
 Suppose that $\D_{\bf i}$ is regular on $A$ and generates $I$.  Then
every row of ${\bf X}$ is a linear combination of rows $i_1,\dots,i_p$
modulo $I_{p+1}({\bf X})$.
 \pf
 Sylvester's relation \Cn1)(2) yields
  $$\D_{\bf i}d_{\bf j}^{\bf l}=\D_{\bf j}d_{\bf i}^{\bf l}.\Dno2.1$$
 Varying ${\bf i}$ in \Cn2.1) yields (1).  On the other hand,
varying ${\bf l}$ in \Cn2.1) yields $\D_{\bf i}I_p({\bf X}_{\bf
j})\subseteq\D_{\bf j}J$.  If $I_p({\bf X}_{\bf j})A=J$, then
$IJ\subseteq\D_{\bf j}J$;  hence, $IJ=\D_{\bf j}J$ because $\D_{\bf
j}\in I$.  Thus (2) holds.
 Finally, \Cn1)(1) yields
  $$\D_{\bf i}{\bf R}_i=\sum\nolimits_{j=1}^p(-1)^{j+p}\D_{{\bf i}i_j}{\bf
R}_{i_j}.$$
 By hypothesis, $\D_{\bf i}$ is regular on $A$ and divides each
$\D_{{\bf i}i_j}$.  Hence (3) holds.

\newsect	 Strongly perfect ideals

\def\seq#1#2{#1_1,\dots,#1_{#2}}

\dfn1
 Let $B$ be a Noetherian ring, $A$ a factor ring of $B$, and $I$ an
$A$-ideal of grade~$g$ such that $\grade_BA/I=s$.  Call $I$ {\it strongly
perfect\/} over $B$ if there exists a generating set $\seq fn$ of $I$ such
that, for $0\le i\le n-g$, the Koszul homology modules $H_i(\seq fn;A)$
are perfect $B$-modules of grade~$s$.

\rmk2
 The notion of strong perfection generalizes Huneke's notion of strong
Cohen--Macaulayness \cite{H2, p.~739}.  Indeed, let $I$ be an ideal of
a local Cohen--Macaulay ring $A$, and write $\widehat A$ as a factor
ring of a regular local ring $B$.  Then $\widehat I$ is strongly
perfect over $B$ if and only if $I$ is strongly Cohen--Macaulay.

Some general results about strong perfection will now be proved.  The
corresponding results about strong Cohen--Macaulayness were proved by
Huneke in \cite{H1} and \cite{H2}.

\lem3
 Let $B$ be a Noetherian ring, $A$ a factor ring of $B$, and $I$ an
$A$-ideal.  Set $s:=\grade_BA/I$  Let $\seq fn$ be an
arbitrary generating set of $I$.
 \part 1 The ideal $I$ is strongly perfect over $B$ if and only if, for
every $i$, the flat dimension over $B$ of $H_i(\seq fn;A)$ is at most
$s$.
 \part 2 If $I$ is strongly perfect over $B$, then the condition in
\Cs1) is satisfied for $\seq fn$.
 \pf
 To prove (1), recall that $IH_i(\seq fn;A)=0$ for all $i$ and that
$H_i(\seq fn;A)\ne0$ if and only if $0\le i\le n-g$ where $g:=\grade
I$.  Hence, $I$ is strongly perfect if the flat dimension of all the
Koszul homology is at most $s$.  Moreover, the converse holds if the
condition in \Cs1) is satisfied for $\seq fn$; so the full converse
follows from (2).

 To prove (2), it suffices to compare a generating set $\seq fn$ with
one of the form $\seq fn,f$.  However, there is a natural isomorphism,
	$$H_i(\seq fn,f;A)=H_i(\seq fn;A)\oplus H_{i-1}(\seq fn;A),$$
 and the assertion follows from the portion of (1)  already proved.

\lem4
 Let $B$ be a Noetherian ring, $A$ a factor ring of $B$, and $I$ an
$A$-ideal.  Let $\seq fn$ be a generating set of $I$.
 \part 1 Let $\seq\D m$ be an $A$-regular sequence contained in $I$,
and let `$\phantom{s}'$' indicate the image in $A':=A/(\seq\D m)$.
Then $I'$ is strongly perfect over $B$ if and only if  $I$  is so.
 \part 2 Let $\seq a r$ be a sequence of elements in $B$ that is
regular on $B$, on $A$, and on $A/I$.  Set $\?B:=B/(\seq a r)B$ and
$\?A:=A/(\seq a r)A$.  Let `$\,\?{\phantom{s}}$' indicate the image in
$\?B$ and in $\?A$.  If $I$ strongly perfect over $B$, then there are
natural isomorphisms,
	$$H_i(\seq {\?f}n;A) = H_i(\seq fn;A)\ox_A\?A,$$
 and $\?I$ is strongly perfect over $\?B$.
 \pf
 To prove (1), we may assume that $m=1$.  Set $\D:=\D_1$.
Because $\D H_j(\seq fn;A)$ vanishes, the exact sequence,
	$$0\To A\TO\D A\To A'\To0,$$
 induces exact sequences
 $$
 0\to H_i(\seq fn;A)\to H_i(\seq{f'}n;A')\to H_{i-1}(\seq fn;A)\to0.$$
 The assertion now follows by induction on $i$ from Lemma~\Cs3)(1).

To prove (2), we may assume that $r=1$.  Set $s:=\grade_BA/I$.  Let $\p$
be an associated prime of the $B$-module $H_i(\seq fn;A)$.  Since $I$ is
strongly perfect over $B$, it follows that $\depth B_{\p}=s$; hence, since
$\p$ is in the support of $A/I$, it follows that $\p$ is associated to
$A/I$; for both these conclusions, see \cite{BVSpgr, (16.17), p.~209}
for example.  Set $a:=a_1$.  Then, therefore, $a$ is regular on
$H_i(\seq fn;A)$.  Now, the exact sequence,
	$$0\To A\TO a A\To \?A\To0,$$
 induces exact sequences,
 $$
 0\to H_i(\seq fn;A)\TO a H_i(\seq{f}n;A)\to H_i(\seq{\?f}n;\?A)\to0.$$
 Hence, they yield the asserted natural isomorphisms.  Furthermore,
$\?I$ is strongly perfect over $\?B$ because $\grade_{\,\?B}\?A/\?I\ge
s$.

\prp5
 Let $B$ be a Noetherian ring, $A$ a factor ring of $B$, and $I$ an
$A$-ideal that is strongly perfect over $B$.  Let $\seq\D
m$ be an $A$-regular sequence contained in $I$, and let $\seq a r$ be a
sequence of elements in $B$ that is regular on $B$, on $A$, and on
$A/I$.  Let `$\,\?{\phantom{s}}$' denote images in  $\?A:=A/(\seq a
r)A$, and assume that $\seq{\?\D}m$ form an $\?A$-regular sequence.
Then, in $\?A$,
   $$\?{\strut(\seq{\D}m)}:\?{\strut I} = \?{\strut(\seq{\D}m):I}.$$
 \pf
 It suffices to verify the asserted equality locally at every
associated prime ideal of the ideal on the right; so we may assume that
all the rings in question are local.
Then, since $\seq a r,\seq\Delta m$ form an $A$-regular sequence,
$\seq\Delta m,\seq a r$ do as well, and hence
 the sequence $\seq a r$ is regular on $A/(\seq{\D}m)$.  Hence, using
Lemma\Cs4)(1), we may reduce to the case $m=0$.  Now, let $\seq fn$ be
a generating set of $I$ with $n\ge1$.  It follows from the definition
of the Koszul complex that there are natural identifications,
	$$0:I=H_n(\seq fn;A) \and \?0:\?I=H_n(\seq{\?f}n;\?A).$$
 Hence the assertion follows from the first assertion of Lemma\Cs4)(2).

\lem6
 Let $B$ be a Noetherian ring, $A$ a factor ring of $B$, and $I$ an
$A$-ideal that is strongly perfect over $B$ with $\grade_BA/I=s$.  Set
$J=0:I$.  Assume that $I+J\not=A$, that $\grade(I+J)\ge1$, and that
$\grade_BA/(I+J)\ge s+1$.  Finally, let `$\,\?{\phantom{s}}$' denote
images in $\?A:=A/J$.  Then $\?I$ is a strongly perfect over $B$ with
$\grade_B\?A/\?I= s+1$.
 \pf Obviously, $J\ne0$.  Hence $\grade I=0$ because $J=0:I$.
Therefore, $I\cap J=0$ because $\grade(I+J)\ge1$.  Let $\seq fn$ be a
generating set of $I$.  Then, by \cite{H2, 1.4, p.~744}, for each $i$,
there is an exact sequence,
   $$0\To\bigoplus J\To H_i(\seq{f}n;A)\To H_i(\seq{\?f}n;\?A)\To0,$$
 where the first term is a direct sum of copies of $J$.  Now, $J=0:I$;
so $J=H_n(\seq{f}n;A)$.  And, the $H_i(\seq{f}n;A)$ have flat dimension
at most $s$ by \Cs3)(1).  Hence the $H_i(\seq{\?f}n;\?A)$ have flat
dimension at most $s+1$.  Hence \Cs3)(1) yields the assertion.

\prp7
Let $R$ be a Noetherian ring.  Let $\bf X$ be a  $p+1$ by $n$ matrix of
variables with $n\ge p+1\ge2$, let $\bf Y$ be the $p+1$ by $p$ matrix
consisting of the first $p$ columns of $\bf X$, and set
	$$B:=R[{\bf X}],\ A:=B/I_{p+1}({\bf X}),\ I:=I_p({\bf Y})A$$
 where $I_{p+1}({\bf X})$ and  $I_p({\bf Y})$ are the ideals of minors
of the indicated sizes.  Then $I$ is an $A$-ideal of grade $1$ that is
strongly perfect over $B$.
 \pf
 Induct on $n$.  Suppose $n=p+1$.  Then $I=I_p({\bf Y})/I_{p+1}({\bf
X})$ where $I_{p+1}({\bf X})$ is generated by a single $B$-regular
element.  On the other hand, Avramov and Herzog \cite{AH, (2.1)(a), p.~252}
proved that $I_p({\bf Y})$ is a strongly perfect $B$-ideal of grade 2.
Hence,  \Cs4)(1) implies that $I$ is an $A$-ideal of grade 1 that
is strongly perfect over $B$.

Suppose  $n\ge p+2$ and that the assertion holds for $n-1$.  Let ${\bf
X}'$  be the matrix consisting of the first $n-1$ columns of ${\bf X}$, set
 $$A':=B/I_{p+1}({\bf X}')B,\ I':=I_p({\bf Y})A',\ J':=I_{p+1}({\bf X})A',$$
 and let $\D'\in A'$ be the image of the $p+1$ by $p+1$ minor of ${\bf
X}$ made of columns $1,\dots,p,n$.  Then $I'$ is an $A'$-ideal of grade
1 that is strongly perfect over $B$ by induction, because the
properties in question are stable under the flat base extension from
$R[{\bf X}']$ to $B$.   Since, moreover, $A'$ is a perfect $B$-module
of grade $n-p-1$, it follows (from  \cite{BVSpgr, (16.18), p.~209} for
example) that
	$$s:=\grade_BA'/I'=\grade_BA'+\grade_{A'}I' =n-p.$$

First, we show that $\D'$ is $A'$-regular.  To this end, let $\q'$  be
an associated prime of $A'$, and let $\q$  be the trace of $\q'$ in $R$.
Since $A'$ is  $R$-flat and $A'/\q A'$ is a domain by
\cite{BVSpgr, (2.10), p.~14}, it follows that $\q'=\q A'$.  Therefore,
$\D'\notin \q'$.

Next, we verify that $\grade_BA'/(I'+J')\ge s+1$ and $\grade
I'+J'\ge2$.  Suppose that the grade of $I_p({\bf
Y})+I_{p+1}({\bf X})$ were equal to that of $I_{p+1}({\bf X})$, which
is $n-p$.  Since the grade of an ideal is the minimum of $\depth B_{\q}$
as $\q$ ranges over all primes containing the ideal, there would be some
$\q$ containing $I_p({\bf
Y})+I_{p+1}({\bf X})$ with $\depth B_{\q}=n-p$.  Since $\q$ also contains
$I_{p+1}({\bf X})$, and that ideal is perfect of grade $n-p$, it
follows that $\q$  would be an associated prime of $I_{p+1}({\bf X})$.
However, an argument like the one above shows that the $B$-ideal
$I_p({\bf Y})$ is not contained in any associated prime of the
$B$-ideal $I_{p+1}({\bf X})$.  Thus
	$$\grade(I_p({\bf Y})+I_{p+1}({\bf X}))
	> \grade I_{p+1}({\bf X}) = n-p.$$
 Therefore, $\grade_BA'/(I'+J')\ge n-p+1=s+1$.  Furthermore, since $A'$
is perfect over $B$ of grade $n-p-1$, it follows (from \cite{BVSpgr,
(16.18), p.~209} for example) that
	$$\grade I'+J'\ge\grade_BA'/(I'+J')-\grade_BA'\ge2.$$

We also have $I'J'\subset(\Delta')$; see \cite{H2, proof of 4.1,
p.~754}.  Indeed, let \hbox{$d_1$, \dots, $d_{p+1}$} denote the
maximal minors of ${\bf Y}$, with alternating signs.  Then in $A'/(\D')$,
	$$(\seq d{p+1}){\bf X}=0,$$
 and hence $I_{p+1}({\bf X})$ annihilates each of the $d_i$ in
$A'/(\D')$.  Therefore, $J'\subseteq(\Delta'):I'$, and equality will
hold if it holds locally at every associated prime $\p$ of the
$B$-module $A'/J'$.  However, since $A'/J'$ is a perfect $B$-module,
$\depth B_{\p}$ is equal to $\grade A'/J'$ (by \cite{BVSpgr, (16.17),
p.~209} for example).  Hence $I'_{\p}=A'_{\p}$ because
	$$\grade A'/J'  <  \grade A'/(I'+J').$$
 Therefore, $J'  =  (\Delta'):I'$ holds locally at $\p$, so globally.

Note that
	$I=(I'+J')/J'\subset A=A'/J'$.
 Factoring out $(\D')$ and using  \Cs4)(1) and \Cs6), we now
conclude that $I$ is strongly perfect over $B$ with
	$$\grade_BA/I=s+1=n-p+1.$$
 But then, by \cite{BVSpgr, (16.18), p.~209} for example,
	$$\grade I=\grade_BA/I-\grade_BA=1.$$

\lem8
 Let $B$ be a Noetherian ring, $A$ a factor ring of $B$, and $I$ an
$A$-ideal of grade $1$ that is strongly perfect over $B$.  Assume that
$A$ is perfect over $B$, and let $J$ be a proper $A$-ideal such that $J\cong
I$.  Then $J$ is an $A$-ideal of grade $1$ that is strongly perfect
over $B$.
 \pf
 Set $s:=\grade_BA/I$.  Since $A$ and $A/I$ are perfect $B$-modules
with $\grade I=1$, it follows (from \cite{BVSpgr, (16.18), p.~209} for
example) that $\grade_BA=s-1$; hence $\grade_BA/J\ge s$.  Furthermore,
$aI=bJ$ for some nonzerodivisors $a$ and $b$ in $A$.  Say $I=(\seq fn)$
and $J=(\seq hn)$ with $af_j=bh_j$.

Let $B_i$ and $Z_i$ denote the modules of boundaries and cycles in the
Koszul complex $(K_\bullet,\partial_\bullet)$.  For every $i$, there
is a commutative diagram
 $$\CD
 K_i(\seq fn;A)   @>\partial_i(\seq fn)>>   K_{i-1}(\seq fn;A)\\
	@|					@VV\mu_a V\\
 K_i(\seq{af}n;A) @>\partial_i(\seq{af}n)>> K_{i-1}(\seq{af}n;A)
 \endCD$$
 where $\mu_a$ denotes multiplication by $a$.  Since $\mu_a$ is
injective, this diagram yields an identification,
	$$Z_i(\seq fn;A)= Z_i(\seq{af}n;A).$$
 Hence, the Isomorphism Theorem yields a natural isomorphism,
	$$B_{i-1}(\seq fn;A)=B_{i-1}(\seq{af}n;A).$$
 Thus there are natural isomorphisms (compare \cite{H1, 1.10
pf., p.~1050}):
	$$\eqalign{Z_i(\seq fn;A)&= Z_i(\seq{af}n;A)\cr
	&=Z_i(\seq{bh}n;A)= Z_i(\seq{h}n;A).\cr}$$
 Likewise, $B_i(\seq fn;A)= B_i(\seq{h}n;A)$.

The $B$-module $A$ has flat dimension $s-1$ because it is perfect of
grade $s-1$.  The $B$-module $H_i(\seq fn;A)$ has flat dimension at
most~$s$ by \Cs3)(1).  It follows by induction on $i$ that $Z_i(\seq
fn;A)$ and $B_i(\seq fn;A)$ have flat dimension at most $s-1$ because
their quotient is $H_i(\seq fn;A)$ and because $Z_i(\seq fn;A)$ is a
first syzygy module of $B_{i-1}(\seq fn;A)$.  Therefore, the above
isomorphisms yield that $H_i(\seq hn;A)$ has flat dimension at most
$s$. But $s\le\grade_BA/J$.  Hence $J$ is strongly perfect over
$B$ by \Cs3)(1), and the proof is complete.

\thm9
 Let $R$ be a Noetherian ring.  Let $\bf X$ be an $m$ by $n$ matrix
with $n\ge m\ge2$ and with entries in $R$.  For $1\le i\le m$, let
${\bf X}_i$ denote the submatrix of $\bf X$ consisting of the last $i$
rows.  Fix $p\ge1$, and assume that the ideals of minors satisfy these
conditions:
	$$\grade I_i({\bf X})=n-i+1\and
	I_i({\bf X})=I_i({\bf X}_i)\for i=p,p+1.$$
 Set $A:=R/I_{p+1}({\bf X})$ and $J:=I_{p}({\bf X})A$.  Denote the
image in $J$ of the minor of ${\bf X}$ with rows $i_1<\dots<i_p$ and
columns $k_1<\dots<k_p$ by $d_{\bf i}^{\bf k}$.
 \part1
 Let ${\bf p}$ be the sequence $m-p+1,\dots,m$.  Then there
exists an $A$-regular element $\D$ of the form
	$\D=\sum_{\bf k}a^{\bf k}d_{\bf p}^{\bf k}$.
 \part2
 Given an $A$-regular element $\D$ as in (1), set $\D_{\bf
i}:=\sum_{\bf k} a^{\bf k}d_{\bf i}^{\bf k}$ and let $I$ be the
subideal of $J$ generated by the various $\D_{\bf i}$.  Then $IJ=\D J$
and $J=(\D):I$.
 \pf
 Obviously, $A=R/I_{p+1}({\bf X}_{p+1})$.  Since $I_{p+1}({\bf
X}_{p+1})$ has generic grade, $A$ is a perfect $R$-module; so $A$ is
grade unmixed (by \cite{BVSpgr, (16.17), p.~209} for example).
Moreover, $\grade I_p({\bf X}_p) >\grade I_{p+1}({\bf X})$.  Therefore,
(1) holds.

Consider (2).  Obviously, (4.7)(2) yields $IJ=\D J$.  So
$J\subseteq(\D):I$.  To prove the opposite inclusion, we may replace
${\bf X}$ by ${\bf X}_{p+1}$.  Indeed, $A$ and $J$ are obviously
unchanged.  Let $I'$ be the ideal generated by the $\D_{\bf i}$ with
$m-p\le i_1$, and suppose $J\supseteq(\D):I'$. Now,
$(\D):I'\supseteq(\D):I$ since $I'\subseteq I$.  Hence $J=(\D):I$.
Thus we may assume $p=m-1$.

Since $J\subseteq(\D):I$, equality will hold if it holds locally at
every associated prime $\q$ of $J$.  Therefore, localizing at $\q$, we
may assume that $R$ is local with $(\D)$, $I$, and $J$ contained in the
maximal ideal of $A$.

The equation $J=(\D):I$ will now be proved in the ``generic'' case and
then specialized.  Let $\bf m$ be the maximal ideal of~$R$, let
$\~{\bf X}$ be an $m$~by~$n$ matrix of indeterminates over $R$, and let $\~B$
denote the localization of the polynomial ring $R[\~{\bf X}]$ at the ideal
$({\bf m},
\~{\bf X}-{\bf X})$.  Let `$\~{\phantom{x}}$' indicate the corresponding
objects
defined using $\~B$ and $\~{\bf X}$ instead of $R$ and ${\bf X}$,
except for $\~J$, which will now denote $I_p(\~{\bf X}_p)\~A$.  Let
${\bf a}$ be the $\~B$-regular sequence consisting of the $mn$ entries
of the difference matrix $\~{\bf X}-{\bf X}$.  Then $\~A/({\bf a})$ is
equal to $A$.  Furthermore, since $\~A$ is a perfect $\~B$-module and
since $\grade_RA$ is equal to $\grade_{\~B}\~A$, it follows that ${\bf
a}$ is $\~A$-regular.  Hence ${\bf a},\~\D$ is $\~A$-regular, and so
$\~\D,{\bf a}$ is $\~A$-regular.  In particular, $\~\D$ is
$\~A$-regular.

Obviously, (4.7)(2) yields $\~I\~J\subseteq(\~\D)$.  Hence,
$\~J\subseteq(\~\D):\~I$, and equality will hold if it holds locally at
every associated prime $\~{\q}$ of $\~J$.  The trace $\q$ of $\~{\q}$ is an
associated prime of $R$, and $\~{\q}/(\~J+\q\~A)$ is an associated prime of
$\~A/(\~J+\q\~A)$; indeed, $\~A/\~J$ is equal to $\~B/I_{p}(\~{\bf
X}_p)$ because $I_p(\~{\bf X}_p)$ contains $I_{p+1}(\~{\bf X})$ as
$m=p+1$, and $\~B/I_{p}(\~{\bf X}_p)$ is (well-known to be) flat over
$R$.  Now, $\~A/(\~J+\q\~A)$ is a domain because it is equal to
$\~B/(I_{p}(\~{\bf X}_p)+\q\~B)$ and the latter is a domain because
$R/\q$ is a domain by \cite{BVSpgr, (2.10), p.~14}.  Therefore,
$\~{\q}=\~J+\q\~A$.

Suppose $\~I\subseteq \~{\q}$.  Then $\~\D_{\bf i}\in\~{\q}$ for every
${\bf i}$.  Take ${\bf i}$ to be the sequence $1,\dots,p$, and pass
momentarily modulo the ideal generated by the last row of $\~{\bf X}$.
Then $\~J$ vanishes, whence $\~\q$ is equal to $\q\~A$.  Hence, since
the $a^{\bf k}$ are the coefficients in the definition of $\~\D_{\bf
i}$, they must lie in $\q$.  Returning to the previous setup, conclude
that $\~\D\in \q\~A$.  Now, $\~A/\q\~A$ is equal to
$\~B/(I_{p+1}(\~{\bf X})+\q\~B)$, which is a domain.  So $\q\~A$ is a
prime.  So it is an associated prime because $\q$ is.  Hence $\~\D$ is
a zerodivisor on $\~A$, contrary to the conclusion drawn above.  Now,
$\~I\not\subseteq \~{\q}$; so $(\~\D):\~I=(\~\D)$ locally at $\~{\q}$.
However, $(\~\D)\subseteq\~J$.  Thus $\~J=(\~\D):\~I$.

To prove that this equation specializes, we first prove that the ideal
$\~I$ has grade 1 and is strongly perfect over $\~B$.  Now,
(4.7)(1) yields
	$$\~ d_{\bf p}^{\lower.2ex\hbox{$\ssize\mkern1mu\bf l$}}\~I
	=\~\D_{\bf p}I_p(\~{\bf X}^{\bf l})\~A$$
 for any ${\bf l}$.  This equation yields an isomorphism of
$\~A$-modules between $\~I$ and $I_{p}(\~{\bf X}^{\bf l})\~A$ because
$\~d_{\bf p}^{\lower.2ex\hbox{$\ssize\mkern1mu\bf l$}}$ and $\~\D_{\bf p}$ are
regular on $\~A$. Furthermore,  \Cs7) implies that $I_p(\~{\bf X}^{\bf l})\~A$
is either the unit ideal or else an $\~A$-ideal of grade 1 that is
strongly perfect over $\~B$.  Hence,  \Cs8) yields that $\~I$ is
an $\~A$-ideal of grade 1 that is strongly perfect over $\~B$.

Finally, $\grade_{\~B}\~A/\~I\le\grade_RA/I$ because $\~I$ has grade
1; hence, because $\~A/\~I$ is perfect, ${\bf a}$ is regular on
$\~A/\~I$.  Proposition~\Cs5) now implies that the equation
$\~J=(\~\D):\~I$ specializes to $J=(\D):I$.  Thus (2) is proved.

\rem10 \enspace If we assume in Theorem \Cs9) that $\grade(d_{\bf p}^{\bf
k})=1$ for some ${\bf k}$, then in (2) we can take $\D=d_{\bf p}^{\bf
k}$ and $I=I_{p}({\bf X}^{\bf k})A$ where ${\bf X}^{\bf k}$ is the $m$
by $p$ submatrix of ${\bf X}$ consisting of columns $k_1<\dots<k_p$.
Moreover, the proof becomes slightly shorter.

\rem11 \enspace Lemmas \Cs4)(1) and \Cs8) yield answers to some
unpublished questions asked by Avramov and Huneke.  Let $B$ be a
Noetherian ring, and $A$ a factor ring that is a perfect $B$-module.
The lemmas imply that, given two $A$-ideals in the same even linkage
class, one ideal is strongly perfect over $B$ if and only if the other
is too; in particular, every $B$-ideal in the linkage class of a
complete intersection is strongly perfect over $B$.  Huneke \cite{H1,
Thm.~1.11, p.~1051} proved the corresponding result for strongly
Cohen-Macaulay ideals in a Gorenstein local ring.

To prove the general case, obviously it suffices to prove the following
assertion.  Let $K$ be a proper $A$-ideal, let $\seq xm$ and $\seq ym$ be
$A$-regular sequences contained in $K$, and set
	$$I:=(\seq xm):K \and J:=(\seq ym):K.$$
 Then $I$ is the unit ideal or is strongly perfect over $B$ if
and only if $J$ is one or the other.

Induct on $m$.  If $m=0$, then $I=J$ and the assertion is trivial.
Suppose $m=1$.  If $I=A$, then $K=(x)$, and so $(x^2):K=(x)$.  Now,
$(x)$ is a proper ideal; moreover, it is strongly perfect because $A$
is perfect.  Hence, we may replace $x$ by $x^2$, and so assume that $I$
is proper.  Similarly, we may assume that $J$ is proper.  Now, in the
total quotient ring of $A$, consider the fractional ideal $xK^{-1}$.
It lies in $A$ because $x\in K$.  Hence $xK^{-1}=(x):K$.  Similarly,
$yK^{-1}=(y):K$.  Therefore, $I$ and $J$ are isomorphic.  Consequently,
the assertion follows from \Cs8).

Suppose $m>1$.  Then we can modify $y_1$ modulo $y_2,\dots,y_m$ so that
$\seq x{m-1},y_1$ form a $A$-regular sequence and $y_1$ is still
$A$-regular (see \cite{H1, proof of Theorem~1.11, p.~1051}).  Set
	$$L:=(\seq x{m-1},y_1):K,\ \?A:=A/(x_1),\and A':=A/(y_1),$$
 and let `$\,\?{\phantom{s}}$' indicate the image in $\?A$ and
`$\phantom{s}'$' that in $A'$.  Then $I$ is the unit ideal or is
strongly perfect over $B$ if and only if $\?I$ is so by \Cs4)(1), if
and only if $\?L$ is so by the induction hypothesis, if and only if
$L'$ is so by \Cs4)(1) applied twice.  Now, $L'$ is so if and only if
$J'$ is so by the induction hypothesis because the ideals $(\seq
x{m-1})$ and $(y_2,\dots,y_m)$ are still generated by $A'$-regular
sequences of length $m-1$ although the given generators need not form
$A'$-regular sequences.  Finally, $J'$ is so if and only if $J$ is so,
by \Cs4)(1) again.  Thus $I$ is the unit ideal or is strongly perfect
over $B$ if and only if $J$ is so, as asserted.

It follows that certain powers of certain ideals $I$ of $B$ have finite
projective dimension; more precisely, if $I$ has grade $m$ and is in
the linkage class of a complete intersection, then $I^i$ has projective
dimension at most~$m+i-2$ in the range $1\le i\le k$ provided that, for
every prime $\p$ containing $I$ with $\depth R_{\p}\le m+k-2$, the number
of generators of $I_\p$ is at most $\depth R_\p$.  Indeed, given a
generating set $\seq fn$ of $I$, set
	$$H_j:=H_j(\seq fn;B) \and S_j:=\hbox{Sym}_j(B^n).$$
 Consider the component $\M_i$ of degree~$i$ of the `approximation
complex' of Simis and Vasconcelos  \cite{SV,
p.~351}:
	$$\M_i: 0\to H_i\ox S_0\to\dots\to H_j\ox S_{i-j}\to
	\dots\to H_0\ox S_i\to0.$$
 These component compexes are acyclic in the range $0\le i\le k-1$ by
the acyclicity lemma because the $B$-modules $H_j$ are either zero or
perfect of grade $m$ and because of the assumption on the number of
generators of each $I_{\p}$; see the proof of Theorem~4.2 in \cite{SV,
p.~353}.  Hence, by the proof of Theorem~4.6 in Herzog, Simis, and
Vasconcelos\cite{HSV, p.~105},
	$$H_0(\M_i) = I^i/I^{i+1} \for 0\le i\le k-1.$$
 Hence, $I^i/I^{i+1}$ has projective dimension at most $m+i$ for $0\le
i\le k-1$ again because the $H_j$ are either zero or perfect of grade
$m$.  Therefore, $I^{i+1}$ has projective dimension at most $m+i-1$ for
$0\le i\le k-1$, as asserted.

\newsect References

\references

AH
 L. Avramov and J. Herzog
 \paper The Koszul algebra of a codimension $2$ embedding
 \mathz 175 1980 249--260

SGAVI
 P. Berthelot et al.,
  ``SGA 6: Th\'eorie des intersections et th\'eor\`eme de Riemann--Roch,''
 Lecture Notes in Math.~{\bf225}, Springer-Verlag,1971

BVMZ
 W. Bruns and U. Vetter
 \paper Length formulas for the local cohomology of exterior powers
 \mathz 191 1986 145--158

BVSpgr
 W. Bruns and U. Vetter,
 ``Determinantal Rings,''
 Lecture Notes in Math.~{\bf 1327}, Springer-Verlag, 1988

B-E
 D. Buchsbaum and D. Eisenbud
 \paper What annihilates a module
 \ja 47 1977 231--43

BR
 D. Buchsbaum and D. Rim
 \paper A generalized Koszul complex. II. Depth and multiplicity
 \tams 111 1964 197--224

Cat
F. Catanese
 \paper Commutative algebra methods and equations of regular surfaces
 \bucharest 68--111

Colley
 S. Colley
 \paper Lines having specified contact with projective varieties
 \vancouver 47--70

Eisenbud
 D. Eisenbud
 \paper Homological algebra on a complete intersection, with an
application to group representations
 \tams 260 1980 35--64

EGAI
 A. Grothendieck, with J. Dieudonn\'e,
``El\'ements de G\'eom\'etrie Alg\'ebrique I,'' Springer-Verlag, 1971

EGAIV
 A. Grothendieck, with J.  Dieudonn\'e, ``El\'ements de
G\'eometrie Alg\'e\-bri\-que $\hbox{\rm IV}_4$,'' Publ.  Math.
I.H.E.S. {\bf 24}, 1965

GP
  L. Gruson and C. Peskine
 \paper Courbes de l'espace projectif: vari\'et\'es de s\'ecantes
 \nice 1--31

HRD
 R. Hartshorne, ``Residues and duality,''
  Lecture Notes in Math.~{\bf20}, Springer-Verlag,1966

HSV
 J. Herzog, A. Simis, and W. V. Vasconcelos
 \paper Koszul homology and blowing-up rings
 \inbook Commutative Algebra -- Proceedings of the Trento Conference
 \bookinfo S. Greco and G.Valla, Eds.
 \publ Marcel Dekker, New York \yr 1983 \pages 79-169

H1
 C. Huneke
 \paper Linkage and the Koszul homology of ideals
 \ajm 104 1982 1043--62

H2
 C. Huneke
 \paper Strongly Cohen--Macaulay schemes
 \tams 277 1983 739--63

TJ
 T. Johnsen
 \paper Eight-secant conics for space curves
 \mathz 211 1992 609--626

Katz
 S. Katz
 \paper Iteration of multiple point formulas and applications to conics
 \inbook Algebraic Geometry -- Sundance, 1986 \bookinfo A. Holme, R.
Speiser, Eds.  Lecture Notes in Math. {\bf 1311} \publ Springer
\yr1988 \pages 147--155

KActa
 S.  Kleiman
 \paper Multiple-point Formulas I: Iteration
 \acta 147 1981 13--49

KSitges
 S.  Kleiman
 \paper Multiple-point formulas II: the Hilbert scheme
 \sitgesii 101--138

KLU
  S. Kleiman, J. Lipman, and B. Ulrich
 \paper The source double-point cycle of a finite map of codimension one
 \paperinfo in ``Complex Projective Varieties,'' G. Ellingsrud, C.
Peskine, G. Sacchiero, and S. A. Stromme (eds.)
 \lmslns 179 1992 199--212

KU
S. Kleiman and B. Ulrich
 \paper Gorenstein algebras, symmetric matrices, self-linked ideals,
and symbolic powers
 \paperinfo Preprint alg-geom/9509005\egroup

Kunz
 E. Kunz
 \book K\"ahler differentials
 \bookinfo  Advandced Lectures in Math.,
 Vieweg, Braunschweig, 1986

L65
 J. Lipman
 \paper Free derivation modules on algebraic varieties
 \ajm 87 1965 874--98

L69
 J. Lipman
 \paper On the Jacobian ideal of the module of differentials
 \pams 21 1969 422--26

McCoy
 N. McCoy
 \book Rings and ideals
 \bookinfo Carus Monograph {\bf 8}, Math. Assoc. of America, 1948

MM
  W. Marar and D. Mond
 \paper Multiple point schemes for corank 1 maps
 \jlms {{\rm(2)} 39} 1989 553--567

Mats
 H. Matsumura
 \book Commutative ring theory
 \bookinfo Cambridge studies in advanced math {\bf 8}, 1986

MP
 D. Mond and R. Pellikaan
 \paper Fitting ideals and multiple points of analytic mappings
 \patzcuaro 107--161

Mum
 D. Mumford,
  ``Lectures on curves on an algebraic surface,''
 Annals of Mathematics Studies {\bf 59},
 Princeton University Press, 1966

Rob
 J. Roberts
 \paper Hypersurfaces with nonsingular normalization and their double
loci \ja 53 1978 253--67

PRob
 P. Roberts \paper Le th\'eor\`eme d'intersection \CR 304
1987 177--180

SV
 A. Simis and W. V. Vasconcelos
 \paper On the dimension and integrality of symmetric algebras
 \mathz 177 1981 341--58

Ulrich
 B. Ulrich,
 ``Algebraic properties of the double-point cycle of a finite map,''
 in preparation

ZN
 R. Zaare-Nahandi
  ``Certain structures on the singular loci at $S_1^{q}$-type singularities,''
 Preprint, University of Tehran, Iran, 1992

\endreferences

  \bigskip

 \eightpoint\smc
 Department of Mathematics, 2--278 MIT, Cambridge, MA 02139,
U.S.A.
\smallskip
 Department of Mathematics, Purdue University, West Lafayette,
IN 47907, U.S.A.
\smallskip
Department of Mathematics, Michigan State University, East Lansing,
MI 48824-1027, U.S.A.

\bye